\pdfoutput=1 

\documentclass[11pt]{amsart}

\usepackage{amssymb,amsthm,enumitem,colonequals,tikz-cd,microtype}
\usepackage[normalem]{ulem} 

\usepackage[osf]{Baskervaldx}
\usepackage[baskervaldx]{newtxmath}
\usepackage[cal=boondoxo]{mathalfa} 

\usepackage[top=3.75cm, bottom=3cm, left=3.5cm, right=3.5cm]{geometry}

\usepackage{xcolor}
\colorlet{darkblue}{blue!55!black}
\colorlet{darkcyan}{cyan!50!black}
\colorlet{darkgreen}{green!60!black}

\PassOptionsToPackage{hyphens}{url}
\usepackage{hyperref}
\hypersetup{
    colorlinks=true,
    linkcolor=darkblue,
    urlcolor=darkcyan,
    citecolor=darkgreen,
}

\def\eqref#1{\textcolor{darkblue}{(\ref{#1})}}

\usepackage[nameinlink]{cleveref} 
\Crefformat{section}{#2\S#1#3}
\Crefmultiformat{section}{#2\S\S#1#3}{ and~#2#1#3}{, #2#1#3}{, and~#2#1#3}

\usepackage[pagewise]{lineno}
\let\oldequation\equation
\let\oldendequation\endequation
\renewenvironment{equation}{\linenomathNonumbers\oldequation}{\oldendequation\endlinenomath}
\expandafter\let\expandafter\oldequationstar\csname equation*\endcsname
\expandafter\let\expandafter\oldendequationstar\csname endequation*\endcsname
\renewenvironment{equation*}{\linenomathNonumbers\oldequationstar}{\oldendequationstar\endlinenomath}
\let\oldalign\align
\let\oldendalign\endalign
\renewenvironment{align}{\linenomathNonumbers\oldalign}{\oldendalign\endlinenomath}
\expandafter\let\expandafter\oldalignstar\csname align*\endcsname
\expandafter\let\expandafter\oldendalignstar\csname endalign*\endcsname
\renewenvironment{align*}{\linenomathNonumbers\oldalignstar}{\oldendalignstar\endlinenomath}

\theoremstyle{plain}
\newtheorem{theorem}{Theorem}[section]
\newtheorem{lemma}[theorem]{Lemma}
\newtheorem{corollary}[theorem]{Corollary}
\newtheorem{proposition}[theorem]{Proposition}

\theoremstyle{definition}
\newtheorem{definition}[theorem]{Definition}
\newtheorem{example}[theorem]{Example}
\newtheorem{remark}[theorem]{Remark}

\newtheorem{convention}[theorem]{Convention}

\newtheorem{reminder}[theorem]{Reminder}

\newtheorem*{ack}{Acknowledgments}

\AddToHook{env/conjecture/begin}{\crefalias{theorem}{conjecture}}
\AddToHook{env/lemma/begin}{\crefalias{theorem}{lemma}}
\AddToHook{env/corollary/begin}{\crefalias{theorem}{corollary}}
\AddToHook{env/proposition/begin}{\crefalias{theorem}{proposition}}
\AddToHook{env/definition/begin}{\crefalias{theorem}{definition}}
\AddToHook{env/remark/begin}{\crefalias{theorem}{remark}}
\AddToHook{env/setup/begin}{\crefalias{theorem}{setup}}
\AddToHook{env/example/begin}{\crefalias{theorem}{example}}
\AddToHook{env/conjecture/begin}{\crefalias{theorem}{Conjecture}}
\AddToHook{env/convention/begin}{\crefalias{theorem}{Convention}}
\AddToHook{env/reminder/begin}{\crefalias{theorem}{Reminder}}

\numberwithin{equation}{section}
\numberwithin{theorem}{section}

\title[Classification, nonexistence, tensor $t$-structures]{Classification and nonexistence \\ for $t$-structures on \\ derived categories of schemes}

\author[A. ~Clark]{Alexander Clark}
\address{A. ~Clark,
School of Mathematics and Statistics, 
Peter Hall Building, University of Melbourne,
Victoria 3010, Australia}
\email{apclark@student.unimelb.edu.au}

\author[P.~Lank]{Pat Lank}
\address{P.~Lank,
Dipartimento di Matematica “F. Enriques”, Universit\`{a} degli Studi di Milano, Via Cesare
Saldini 50, 20133 Milano, Italy}
\email{plankmathematics@gmail.com}

\author[K. ~Manali Rahul]{Kabeer Manali Rahul}
\address{K. ~Manali Rahul,
Center for Mathematics and its Applications, 
Mathematical Science Institute, Building 145, 
The Australian National University, 
Canberra, ACT 2601, Australia}
\email{kabeer.manalirahul@anu.edu.au}

\author[C. ~J. Parker]{Chris J. Parker}
\address{C. ~J. Parker,
Fakult\"at f\"ur Mathematik, 
Universit\"at Bielefeld, 33501, Bielefeld, Germany}
\email{cparker@math.uni-bielefeld.de}

\date{\today}



\begin{document}
    
\begin{abstract}
    Given a suitable Noetherian scheme, we classify tensor $t$-structures on the bounded derived category of coherent sheaves and its variants with prescribed support. Furthermore, we show that the existence of such $t$-structures restricting to perfect complexes detects regularity, recovering a theorem of Neeman in the affine case by different methods. Our tools establish local-to-global principles for tensor $t$-structures.
\end{abstract}

\maketitle

\section{Introduction}
\label{sec:introduction}

\subsection{What is known}
\label{sec:introduction_known}

The existence and classification of $t$-structures on triangulated categories arising in algebraic geometry has given new ways to extract geometric information from the homological data. A notable application of $t$-structures is a conjecture of Antieau--Gepner--Heller \cite{Antieau/Gepner/Heller:2019}. It predicts that $\operatorname{Perf}(X)$ for a Noetherian scheme $X$ admits a bounded $t$-structure if, and only if, $X$ is regular. This conjecture has now been resolved: first in the affine case by Smith \cite{Smith:2022}, and in full generality by Neeman \cite{Neeman:2022}.

The classification problem for $t$-structures on triangulated categories associated to schemes is subtle. In this work, we consider categories such as $D^b_{\operatorname{coh}}(X)$, $D_{\operatorname{qc}}(X)$, $\operatorname{Perf}(X)$, and their variants with support in a closed subset. A complete classification of all $t$-structures on $D_{\operatorname{qc}}(X)$ is out of reach even in the affine case \cite{Stanley:2010}. Instead, one is naturally led to restrict attention to well-behaved classes of $t$-structures. Key advances in this direction include the classification of compactly generated $t$-structures on $D_{\operatorname{qc}}(X)$ for affine schemes \cite{AlonsoTarrio/JeremiasLopez/Saorin:2010, Hrbek:2020}. 

Passing to the global (nonaffine) setting introduces additional difficulties. In particular, many techniques from the affine case no longer apply directly. A robust class of $t$-structures in this context are those that are both \textit{compactly generated} and \textit{tensor compatible}. The classification of such $t$-structures on $D_{\operatorname{qc}}(X)$ for the global case has been achieved in \cite{Dubey/Sahoo:2023, Lank:2025}. These results establish a bijective correspondence between compactly generated tensor $t$-structures on $D_{\operatorname{qc}}(X)$ and Thomason filtrations of $X$.

\subsection{What we do}
\label{sec:introduction_what_we_do}

The technical core of the paper concerns the local-to-global behavior of compactly generated tensor $t$-structures on $D_{\operatorname{qc}}(X)$. A central role is played by the classification of such $t$-structures in terms of Thomason filtrations. Our approach shows that one can effectively study $t$-structures under passage to open subschemes, local rings, and to subcategories defined by support conditions. This allows us to reduce global questions to local conditions on stalks and topological conditions. See \Cref{sec:key_ingredients} for details on the tools that we develop.

\subsubsection{Relative weak Cousin}
\label{sec:introduction_classification}

Our first result concerns compactly generated tensor $t$-structures on $D_{\operatorname{qc}}(X)$ that restrict to $D^b_{\operatorname{coh}}(X)$ and its variants with prescribed support. We give a topological characterization, in terms of Thomason filtrations, for when such restrictions occur. In the affine setting, related classification results were established in \cite{AlonsoTarrio/JeremiasLopez/Saorin:2010} for schemes admitting a dualizing complex, and later in \cite{Takahashi:2023} for affine CM-excellent schemes of finite Krull dimension.

We now state our first main result.

\begin{proposition}
    \label{prop:weak_cousin_cm_excellent_stalks}
    Suppose $X$ is a quasi-compact CM-excellent scheme of finite Krull dimension. A Thomason filtration $\phi$ on $X$ is weak Cousin across $Z$ if, and only if, its associated $t$-structure on $D_{\operatorname{qc},Z}(X)$ restricts to $D^b_{\operatorname{coh},Z}(X)$.
\end{proposition}

\Cref{prop:weak_cousin_cm_excellent_stalks} provides a characterization which applies to complexes supported on a fixed closed subset. This feature is absent from the affine results in the references above. A key ingredient in the proof is a relative form of the weak Cousin condition for Thomason filtrations, extending the notion of \cite{AlonsoTarrio/JeremiasLopez/Saorin:2010}. See \Cref{sec:results1}. In the affine absolute case (i.e.\ $Z=X=\operatorname{Spec}(R)$), \cite[Corollary 5.6]{Hrbek/Martini:2024} says that these results and their hypothesis of CM-excellence are optimal. It is natural to ask whether a global version of loc.\ cit.\ holds.

The intuition behind the proof of \Cref{prop:weak_cousin_cm_excellent_stalks} begins with studying how Thomason filtrations on $X$ behave upon restriction to a closed subset. In particular, being relative weak Cousin over $Z$ simply means that the restricted Thomason filtration satisfies the weak Cousin condition when viewed as a Thomason filtration on $Z$. 

Our next result yields a classification of all tensor $t$-structures on $D^b_{\operatorname{coh},Z}$:

\begin{theorem}
    \label{thm:weak_cousin_cm_excellent_stalks}
    Let $X$ be a quasi-compact CM-excellent scheme of finite Krull dimension. For any closed $Z\subseteq X$, there is a one-to-one correspondence:
    \begin{displaymath}
        \begin{aligned}
            \{ \otimes\textrm{-aisle on } & D^b_{\operatorname{coh},Z}(\mathcal{X})\} 
            \iff \{ \textrm{Thomason filtrations on } X \textrm{ that are weak Cousin across } Z \}.
        \end{aligned}
    \end{displaymath}
\end{theorem}

A key ingredient to proving this result is \Cref{thm:weak_cousin_cm_excellent_stalks}. Furthermore, \Cref{thm:weak_cousin_cm_excellent_stalks} is quite applicable, e.g.\ Noetherian schemes admitting dualizing complexes \Cref{ex:dualizing_complex_implies_cm_excellent_finite_krull_dim}. Its ramifications extend to perverse sheaf theory. Particularly, our result can be viewed as an effective conclusion to classifying $t$-structures on $D^b_{\operatorname{coh}}(X)$ obtained from perversity functions on $X$. See \cite{Arikin/Bezrukavnikov:2010}, \cite[Theorem 9.1]{Gabber:2004}, \cite[Theorem 5.9]{Kashiwara:2004}, and \cite{Sahoo:2024}.

\subsubsection{Restriction to perfects}
\label{sec:introduction_nonexistence}

We next study $t$-structures on $D_{\operatorname{qc},Z}(X)$ that restrict to $\operatorname{Perf}_Z(X)$ for a fixed closed subset $Z \subset X$. A key point is that restrictability implies a support theoretic constraint in terms of the regular locus (see \Cref{prop:restrictability_implies_regular_locus_containment}). This leads to the following:

\begin{theorem}
    \label{thm:appearance_for_restrictability_via_filtration}
    A nonempty closed subset $Z$ of a Noetherian scheme $X$ is contained in the regular locus of $X$ if, and only if, there is a Thomason filtration on $X$ whose associated $t$-structure on $D_{\operatorname{qc}}(X)$ restricts to $\operatorname{Perf}_Z (X)$ and is not constant on each connected component $Z_t$ of $Z$.
\end{theorem}

This recovers \cite[Theorem 0.1]{Neeman:2022} in the affine case and extends \cite{Smith:2022} to arbitrary Krull dimension. Using a result of Hochster \cite{Hochster:1973}, it also applies to explicit examples of Noetherian integral domains of infinite Krull dimension that are not regular (see \Cref{ex:nagata_singular}). The proof uses local algebra techniques rather than the methods of \cite{Neeman:2022} involving approximable triangulated categories.

As a consequence, we obtain the following classification of $t$-structures on $\operatorname{Perf}_Z(X)$:

\begin{corollary}
    \label{cor:classification_for_relative_perf}
    Let $X$ be a quasi-compact CM-excellent scheme of finite Krull dimension. Suppose $Z$ is a nonempty closed subset of $X$ with connected components $Z_t$. A Thomason filtration $\phi$ on $X$ has an associated $t$-structure on $D_{\operatorname{qc}}(X)$ which restricts to $\operatorname{Perf}_Z (X)$ if, and only if, $\phi\cap Z_t$ satisfies weak Cousin for $Z_t$ contained in the regular locus of $X$ and $\phi\cap Z_t$ is constant on $Z_t$ otherwise.
\end{corollary}

\begin{ack}
    Lank, Manali Rahul, and Parker were supported under the ERC Advanced Grant 101095900-TriCatApp.\ Manali Rahul was supported under the Australian Research Council Grant DP200102537, is a recipient of an Australian Government Research Training Program Scholarship, and would like to thank Universit\`{a} degli Studi di Milano for their hospitality during his stay there. Additionally, Parker thanks the University of Milan for a visit, and was support by the Deutsche Forschungsgemeinschaft (SFB-TRR 358/1 2023 - 491392403). Lank was partly supported by the National
    Science Foundation under Grant No.\ DMS-2302263. Clark was supported by an Australian Government Research Training Program Scholarship. The authors would like to Michal Hrbek for very useful comments in an earlier version of our work. Lank would like to thank Takumi Murayama for helpful commutative algebra discussions regarding \cite{Hochster:1973, Nagata:1962}.
\end{ack}

\section{Preliminaries}
\label{sec:prelim}

\subsection{Notation/convention}
\label{sec:prelim_notation}

Let $X$ be a Noetherian scheme. Denote by $D_{\operatorname{qc}}(X)$ the derived category of complexes of $\mathcal{O}_X$-modules with quasi-coherent cohomology. Let $D^b_{\operatorname{coh}}(X)$ be the strictly full subcategory of objects in $D_{\operatorname{qc}}(X)$ with coherent and bounded cohomology. Set $\operatorname{Perf}(X)$ to be the strictly full subcategory of $D_{\operatorname{qc}}(X)$ consisting of perfect complexes. If $Z$ is a subset of $X$, then $D_{\operatorname{qc},Z}(X)$ (resp. $D^b_{\operatorname{coh},Z}(X)$, $\operatorname{Perf}_Z (X)$) is the full subcategory in (resp. $D^b_{\operatorname{coh}}(X)$, $\operatorname{Perf} (X)$) consisting of objects whose support is contained in $Z$. In the special case $X=\operatorname{Spec}(R)$, we will abuse notation and write $D_{\operatorname{qc}}(R)$ for $D_{\operatorname{qc}}(X)$; similarly for the other subcategories. If $X=\operatorname{Spec}(R)$ and $I$ is an ideal in $R$, then we will write $K_R(I)$ for the Koszul complex on a set of generators of $I$. That is, $K_R(I)= \otimes^n_{j=1} K_R(r_j)$ where $K_R (r_j)=\operatorname{cone}(R \xrightarrow{r_j \cdot - } R)$ and $I=(r_1,\ldots,r_n)$. Note that $K_R(I)$ is in $\operatorname{Perf}(X)$ for all ideals $I$ in $R$. The $n$-th cohomology of an object $E$ in $D_{\operatorname{qc}}(X)$ is denoted by $\mathcal{H}^n (E)$. Additionally, if $X$ is affine, then we write this instead as $H^n (E)$. There are the following triangulated equivalences of $D_{\operatorname{qc}}(X)$ with $D(\operatorname{Qcoh}(X))$, and $D^b_{\operatorname{coh}}(X)$ with $D^b (\operatorname{coh}(X))$, which respectively appear in \cite[\href{https://stacks.math.columbia.edu/tag/09T4}{Tag 09T4}]{StacksProject} and \cite[\href{https://stacks.math.columbia.edu/tag/0FDB}{Tag 0FDB}]{StacksProject}. 

\subsection{t-structures}
\label{sec:t_structures}

We discuss $t$-structures and related constructions within triangulated categories. See \cite{Beilinson/Berstein/Deligne/Gabber:2018} for further details. Let $\mathcal{T}$ be a triangulated category equipped with a shift functor $[1]\colon \mathcal{T} \to \mathcal{T}$, and assume it admits all small coproducts. Recall a subcategory is `strictly full' if it is a full subcategory which is closed under isomorphisms. If $\mathcal{S}$ is a subcategory of $\mathcal{T}$, then $\operatorname{Coprod}(\mathcal{S})$ is the smallest strictly full subcategory of $\mathcal{T}$ containing $\mathcal{S}$ that is closed under extensions and small coproducts.

A pair of strictly full subcategories $(\mathcal{A},\mathcal{B})$ in $\mathcal{T}$ is called a \textbf{$t$-structure} if the following conditions are satisfied:
\begin{enumerate}
    \item $\operatorname{Hom}(A,B)=0$ for all $A$ in $\mathcal{A}$ and $B$ in $\mathcal{B}[-1]$
    \item $\mathcal{A}, \mathcal{B}$ are respectively closed under positive and negative shifts
    \item For any object $E$ in $\mathcal{T}$, there is a distinguished triangle with $A$ in $\mathcal{A}$ and $B$ in $\mathcal{B}[-1]$:
    \begin{displaymath}
        A \to E \to B \to A[1].
    \end{displaymath}
\end{enumerate}
A \textbf{preaisle} on $\mathcal{T}$ is a strictly full subcategory $\mathcal{A}$ that is closed under positive shifts and extensions; whereas an \textbf{aisle} on $\mathcal{T}$ is a preaisle whose inclusion $\mathcal{A} \to \mathcal{T}$ admits a right adjoint. If the aisle exists, denote by $\overline{ \langle \mathcal{S} \rangle}^{(-\infty,0]}$ the smallest cocomplete aisle of $\mathcal{T}$ containing a subcategory $\mathcal{S}$. A $t$-structure $(\mathcal{A},\mathcal{B})$ on $\mathcal{T}$ is said to \textbf{restrict} to a triangulated subcategory $\mathcal{S}$ of $\mathcal{T}$ if the pair $(\mathcal{A}\cap \mathcal{S},\mathcal{B} \cap \mathcal{S})$ is a $t$-structure on $\mathcal{S}$. In a similar vein, we can speak of aisles restricting if their associated $t$-structure does so.

\begin{lemma}\label{lem:lifitng_tensor_structure_from_perfect_categorie_abstract}
    Let $\mathcal{T}$ be a compactly generated triangulated category which admits all small coproducts. Denote the triangulated subcategory of compact objects by $\mathcal{T}^c$. Then any $t$-structure on $\mathcal{T}^c$ is the restriction of a compactly generated $t$-structure on $\mathcal{T}$.
\end{lemma}

\begin{proof}
    This is essentially \cite[Lemma 1.8 \& Proposition 1.9]{Neeman:2021} coupled with the fact that aisles are closed under direct summands (see \cite[Corollary 1.4]{AlonsoTarrio/JeremiasLopez/Salorio/Souto:2003}). Indeed, for any aisle $\mathcal{A}$ on $\mathcal{T}^c$, one can lift to $\operatorname{Coprod}(\mathcal{A})$ which is a compactly generated aisle on $\mathcal{T}$ (see \cite[Theorem 2.6]{Canonaco/Neeman/Stellari:2025} for details). In fact, we have that $\overline{ \langle \mathcal{A} \rangle}^{(-\infty,0]} = \operatorname{Coprod}(\mathcal{A})$.
\end{proof}

The distinguished triangles in the datum of a $t$-structure are called the \textbf{truncation triangles}, and are unique up to unique isomorphism of triangles. If $A\to E \to B \to A[1]$ is a truncation triangle, then we typically define $\tau^{\leq 0} E:=A$ and $\tau^{\geq 1}:=B$ respectively as the `connective' and `coconnective' components of $E$. In fact, these assignments give us well-defined functors $\tau^{\leq 0} \colon \mathcal{T}\to \mathcal{A}$ and $\tau^{\geq 1} \colon \mathcal{T} \to \mathcal{B}$, which are respectively called the \textbf{connective} and \textbf{coconnective truncation functors}. Given a $t$-structure $(\mathcal{A},\mathcal{B})$, the subcategory $\mathcal{A}$ is an aisle. On the other hand, for an aisle $\mathcal{A}$, there is a $t$-structure $(\mathcal{A},\mathcal{B})$ where $\mathcal{B} \colonequals \mathcal{A}^{\perp}[1]$, see \cite{Keller/Vossieck:1988}. A $t$-structure may also be denoted by $(\mathcal{T}^{\leq 0}, \mathcal{T}^{\geq 0})$, and occasionally we use this convention. If $(\mathcal{T}^{\leq 0}, \mathcal{T}^{\geq 0})$ is a $t$-structure on $\mathcal{T}$ and $n$ is an integer, then the pair $(\mathcal{T}^{\leq n}, \mathcal{T}^{\geq n})$ is also a $t$-structure on $\mathcal{T}$ where $\mathcal{T}^{\leq n}:= \mathcal{T}^{\leq 0}[-n]$ and $\mathcal{T}^{\geq n}:= \mathcal{T}^{\geq 0}[-n]$. 

\begin{example}\label{ex:standard_t_structure}
    Let $\mathcal{A}$ be an abelian category. Consider the following categories: $D^{\leq 0}(\mathcal{A})$ consists of objects $E$ in $D(\mathcal{A})$ such that strictly positive cohomology vanishes, i.e.\ $H^j (E)$ is zero for $j >0$, and $D^{\geq 0}(\mathcal{A})$ consists of objects $E$ in $D(\mathcal{A})$ such that strictly negative cohomology vanishes, i.e.\ $H^j (E)$ is zero for $j <0$. The pair $(D^{\leq 0}(\mathcal{A}), D^{\geq 0}(\mathcal{A}))$ is called the \textbf{standard $t$-structure} on $D(\mathcal{A})$. An example of interest in our work is $\mathcal{A}=\operatorname{Qcoh}(X)$.
\end{example}

\begin{example}
    Let $R$ be a commutative Noetherian ring and $Z$ be a Thomason subset of $\operatorname{Spec}(R)$. Denote by $\mathbf{R}\Gamma_Z \colon D_{\operatorname{qc}}(R)\to D_{\operatorname{qc}}(R)$ the \textbf{local cohomology} functor associated to $Z$. This is the connective truncation functor associated to the $t$-structure arising from the constant Thomason filtration at $Z$. See \cite[Example 3.4]{Smith:2022}. The $i$-th cohomology of an object $E$ with respect to this $t$-structure is denoted by $H^i_Z (E)$.
\end{example}

An aisle $\mathcal{U}$ on $\mathcal{T}$ is \textbf{compactly generated} if there exists a collection of compact objects $\mathcal{P}$ in $\mathcal{T}$ such that $\overline{ \langle \mathcal{P} \rangle}^{(-\infty,0]} = \mathcal{U}$, and a $t$-structure on $\mathcal{T}$ is \textbf{compactly generated} if its associated aisle is such. Assume $(\mathcal{T},\otimes,1)$ is a tensor triangulated category. If $\mathcal{S}$ is a preaisle on $\mathcal{T}$ such that $\mathcal{S}\otimes \mathcal{S}$ is contained in $\mathcal{S}$ and $1$ is an object of $\mathcal{S}$, then a \textbf{tensor aisle}, also known as a \textbf{$\otimes$-aisle}, with respect to $\mathcal{S}$ is an aisle $\mathcal{U}$ such that $\mathcal{S} \otimes \mathcal{U}$ is contained in $\mathcal{U}$. Moreover, a $t$-structure on $\mathcal{T}$ is said to be a \textbf{tensor $t$-structure} with respect to $\mathcal{S}$ if its associated aisle is a $\otimes$-aisle with respect to $\mathcal{S}$. If it is clear from context, we will drop the terminology `with respect to $\mathcal{S}$'. For a given subcategory $\mathcal{C}$ of $\mathcal{T}$, the smallest $\otimes$-aisle containing $\mathcal{C}$ in $\mathcal{T}$ is denoted $\overline{ \langle \mathcal{C} \rangle}^{(-\infty,0]}_{\otimes}$. In fact, since $\mathcal{T}$ is compactly generated, \cite[Proposition 3.3]{Hrbek/Lank/LeGros/Pavon:2026} ensures that $\overline{ \langle \mathcal{C} \rangle}^{(-\infty,0]}_{\otimes}$ is an aisle.

\begin{convention}
    Let $Z$ be a closed subset of a Noetherian scheme $X$. A $t$-structure $(\mathcal{A},\mathcal{B})$ on $\operatorname{Perf}_Z (X)$ (resp.\ $D^b_{\operatorname{coh},Z}(X)$, $D_{\operatorname{qc},Z}(X)$) is said to be \textbf{tensor} if $(\operatorname{Perf}_Z (X) \cap D^{\leq 0}_{\operatorname{qc}}(X) )\otimes^{\mathbf{L}} \mathcal{A} \subseteq \mathcal{A}$ (resp.\ $D^{\leq 0}_{\operatorname{coh}}(X) \otimes^{\mathbf{L}}\mathcal{A} \subseteq \mathcal{A}$, $D^{\leq 0}_{\operatorname{qc},Z}(X) \otimes^{\mathbf{L}} \mathcal{A} \subseteq \mathcal{A}$). 
\end{convention}

\subsection{Thomason filtrations}
\label{sec:thomason_ftilration}

Let $X$ be a Noetherian scheme. We briefly discuss Thomason filtrations on such spaces. Recall that \textbf{Thomason subset} of $X$ is an arbitrary union of subsets that are both closed and have a quasi-compact complement. The fact $X$ is Noetherian ensures that this condition is equivalent to being a specialization  closed subset of $X$. Let $\operatorname{Spcl}(X)$ be the collection of specialization  closed subsets of $Y$. A function $\phi \colon \mathbb{Z} \to \operatorname{Spcl}(X)$ is called a \textbf{Thomason filtration} on $X$ if for each integer $n$ we have that $\phi(n)$ is a Thomason subset of $X$, and $\phi(n+1)$ is contained in $\phi (n)$. It is possible to `shift' a Thomason filtration $\phi$ by an integer $s$ as follows: $\phi[s]$ is determined by assignment $n\mapsto \phi (n+s)$. 


\begin{example}\label{ex:truncate_thomason_filtration}
    For each Thomason filtration $\phi$ on $X$ and for each integer $k$ one can associate another Thomason filtration $\phi^{\geq k}$ which is constant at $\phi (k)$ for $i\leq k$ as follows:
    \begin{displaymath}
        \phi^{\geq k}(i):=\begin{cases}\phi (k),&\text{for }i\leq k,\\ \phi (i),&\text{for } i\geq k\end{cases}
    \end{displaymath}
\end{example}

\begin{definition}
    Let $S$ be a subset of $Y$ and $\phi$ be a Thomason filtration on $Y$.
    \begin{enumerate}
        \item $\phi$ is \textbf{eventually vanishing} if the intersection of all $\phi(n)$ is empty.
        \item The \textbf{restriction} of $\phi$ to $S$, denoted $\phi\cap S$, is given by $n\mapsto \phi (n)\cap S$.
        \item $\phi$ is \textbf{constant} on $S$ if $( \phi \cap S)(n)$ is $S$ for all $n$ or $\emptyset$ for all $n$.
    \end{enumerate}
\end{definition}

\begin{example}\label{ex:standard_t_structure_filtration}
    Let $Z$ be a closed subset of $X$. The standard $t$-structure on $D^b_{\operatorname{coh},Z}(X)$ is determined by the Thomason filtration on $X$ given $n\mapsto Z$ for $n\leq 0$ and $n\mapsto \emptyset$ for $n>0$. This is a Thomason filtration which is nonconstant over $Z$. We denote this by $\phi_{\operatorname{st}}$ in the special case $Z=X$.
\end{example}

\begin{definition}\label{def:weak_cousin}
    \hfill
    \begin{enumerate}
        \item A generalization  $p\rightsquigarrow q$ (i.e.\ $q\in \overline{\{ p\}}$) is \textbf{direct} if $p\not=q$ and there does not exist a point $s$ in $Y$ such that $p\rightsquigarrow s \rightsquigarrow q$ where $p\not=s$.
        \item A Thomason filtration $\phi \colon \mathbb{Z} \to \operatorname{Spcl}(Y)$ is said to satisfy the \textbf{weak Cousin condition} if for every integer $n$, the set $\phi (n-1)$ contains all direct generalization s of each point in $\phi (n)$ (i.e.\ for each direct generalization  $p \rightsquigarrow q$ with $q$ in $\phi (n)$, one has $p$ is in $\phi (n-1)$).
    \end{enumerate}
\end{definition}

\begin{reminder}
    \label{rem:classification}
    There is a bijective correspondence between the collection of compactly generated tensor $t$-structures on $D_{\operatorname{qc}}(X)$ and the collection of Thomason filtrations on $X$. Initially, this was proven for Noetherian affine schemes \cite[Theorem 3.11]{AlonsoTarrio/JeremiasLopez/Saorin:2010}, and later extended to Noetherian schemes \cite[Theorem 5.11]{Dubey/Sahoo:2023}. However, it has been fully generalized to affine schemes in \cite[Theorem 5.6]{Hrbek:2020}, and recently extended to algebraic stacks \cite[Theorem 1.1]{Lank:2025} using loc.\ cit. with methods independent than that of \cite{Dubey/Sahoo:2023}. We briefly recall the assignments.
    
    Let $\phi$ be a Thomason filtration on $X$. The associated category $\mathcal{U}_\phi$ is a compactly generated $\otimes$-aisle on $D_{\operatorname{qc}}(X)$, for which we will denote the corresponding compactly generated tensor $t$-structure by $(\mathcal{U}_\phi, \mathcal{V}_\phi)$ and the truncation functors by $\tau_\phi ^{\leq 0}$ and $ \tau_\phi ^{\geq 1}$. Specifically, $\mathcal{U}_\phi$ is the strictly full subcategory of $D_{\operatorname{qc}}(X)$ consisting of objects $E$ such that $\operatorname{supp}(\mathcal{H}^i(E))$ is contained in $\phi (i)$ for all integers $i$. The associated aisle $\mathcal{U}_\phi$ is compactly generated by the following subcategory:
    \begin{displaymath}
        \bigcup_{n \in \mathbb{Z}} \big(\operatorname{Perf}(X) \cap D_{\operatorname{qc},\phi (n)}^{\leq n}(X)\big).
    \end{displaymath}
    
    Let $\mathcal{U}$ be a $\otimes$-aisle on $D_{\operatorname{qc}}(X)$ that is compactly generated by a collection $\mathcal{S}$. The associated Thomason filtration, denoted by $\phi_{\mathcal{U}}$, is given by the following:
    \begin{displaymath}
        \phi_{\mathcal{U}}(i)= \bigcup_{P \in \mathcal{S}} \bigcup_{j\geq i} \operatorname{supp}(\mathcal{H}^j (P)).
    \end{displaymath}
    
    We may impose a partial ordering on $\operatorname{Spcl}(X)$, which naturally induces a partial ordering on the collection of Thomason filtrations on $X$. Moreover, there is a partial ordering on the collection of compactly generated tensor $t$-structures on $D_{\operatorname{qc}}(X)$ by inclusion on the aisles. With respect to these orders, the bijective correspondence is order preserving.
\end{reminder}

\begin{convention}
    Let $\phi$ be a Thomason filtration on $X$. We at times call the corresponding compactly generated tensor $t$-structure its \textbf{associated $t$-structure on $D_{\operatorname{qc}}(X)$}, and vice-versa for the $t$-structure corresponding to $\phi$. These notions extend in a similar vein for aisles associated to Thomason filtrations.
\end{convention}

\begin{lemma}
    \label{lem:dbcoh_tensor_via_perfect_tensor}
    Let $X$ be a Noetherian scheme. Then an aisle $\mathcal{A}$ on $D^b_{\operatorname{coh}}(X)$ being tensor is equivalent to it being closed under the tensor action by $\operatorname{Perf}(X)\cap D^{\leq 0}_{\operatorname{qc}}(X)$. 
\end{lemma}

\begin{proof}
    Since $\operatorname{Perf}(X)\cap D^{\leq 0}_{\operatorname{qc}}(X)\subseteq D^{\leq 0}_{\operatorname{coh}}(X)$, it follows that $\mathcal{A}$ being tensor aisle on $D^b_{\operatorname{coh}}(X)$ implies
    \begin{displaymath}
        (\operatorname{Perf}(X)\cap D^{\leq 0}_{\operatorname{qc}}(X) ) \otimes^{\mathbf{L}} \mathcal{A} \subseteq \mathcal{A}.
    \end{displaymath}
    To see the converse, let $E\in D^{\leq 0}_{\operatorname{coh}}(X) \otimes^{\mathbf{L}} \mathcal{A}$ where $\mathcal{A}$ is an aisle on $D^b_{\operatorname{coh}}(X)$ which is closed under the tensor action by $\operatorname{Perf}(X) \cap D^{\leq 0}_{\operatorname{qc}}(X)$. Consider the truncation triangle of $E$ with respect to $\mathcal{A}$ on $D^b_{\operatorname{coh}}(X)$,
    \begin{displaymath}
        \tau^{\leq 0} E \to E \to \tau^{\geq 1} E \to \tau^{\leq 0} E.
    \end{displaymath}
    By \cite[Corollary 3.4 \& Theorem 3.10]{Hrbek/Lank/LeGros/Pavon:2026}, $\overline{\langle \mathcal{A} \rangle}^{(-\infty,0]}_{\otimes}$ exists and $(\operatorname{Perf}(X) \cap D^{\leq 0}_{\operatorname{qc}}(X))\otimes^\mathbf{L} \mathcal{A}$ compactly generates it as an aisle on $D_{\operatorname{qc}}(X)$. Now, the distinguished triangle above corresponds to the truncation triangle for $E$ with respect to $\overline{\langle \mathcal{A} \rangle}^{(-\infty,0]}_{\otimes}$ on $D_{\operatorname{qc}}(X)$. Indeed, $\tau^{\geq 1} E \in \mathcal{A}^\perp \cap D^b_{\operatorname{coh}}(X)$, and so $E\in (\overline{\langle \mathcal{A} \rangle}^{(-\infty,0]}_{\otimes})^\perp \cap D^b_{\operatorname{coh}}(X)$ (e.g.\ use \cite[Lemma 3.1]{AlonsoTarrio/JeremiasLopez/Salorio/Souto:2003}). Since $E\in D^{\leq 0}_{\operatorname{coh}}(X) \otimes^{\mathbf{L}} \mathcal{A}$, it follows that $E\in \overline{\langle \mathcal{A} \rangle}^{(-\infty,0]}_{\otimes}$. Thus, the morphism  $\tau^{\leq 0} E \to E$ is an isomorphism, which implies $E\in \mathcal{A}$ as desired.
\end{proof}

\section{Key ingredients}
\label{sec:key_ingredients}

This section examines the local-to-global behaviour of compactly generated $t$-structures on $D_{\operatorname{qc}}(X)$ where $X$ is a fixed Noetherian scheme. We begin by establishing a few key ingredients before proceeding to proving our main results.

\begin{lemma}\label{lemma:AJS2010_4.2_generalization }
    Suppose $\phi$ is a Thomason filtration on $X$. Let $j$ be an integer. Then $\tau^{\leq 0}_\phi D^{\geq j}_{\operatorname{qc}}(X)$ and $\tau^{\geq 1}_\phi D^{\geq j}_{\operatorname{qc}}(X)$ are contained in $D^{\geq j}_{\operatorname{qc}}(X)$.
\end{lemma}

\begin{proof}
    This is essentially the proof of \cite[Lemma 4.2]{AlonsoTarrio/JeremiasLopez/Saorin:2010} where in the global setting one invokes \cite[Theorem 3.6]{Lank:2025} and \cite[Example 2.1]{Hrbek/Lank/LeGros/Pavon:2026}.
\end{proof}

\begin{proposition}\label{prop:filtrations_agreeing_after_some_point}
    Let $\phi$ be a Thomason filtration on $X$. Suppose $F$ is an object of $D^{\geq k}_{\operatorname{qc}}(X)$. Then $\tau^{\leq 0}_\phi F$ is isomorphic to $\tau^{\leq 0}_\psi F$ in $D_{\operatorname{qc}}(X)$ where $\psi:=\phi^{\geq k}$ (see \Cref{ex:truncate_thomason_filtration}).
\end{proposition}

\begin{proof}
    Let $F$ be an object of $D^{\geq k}_{\operatorname{qc}}(X)$. It suffices to check the truncation triangle
    \begin{displaymath}
        \tau^{\leq 0}_\phi F\to F\to \tau^{\geq 1}_\phi F \to (\tau^{\leq 0}_\phi F)[1]
    \end{displaymath}
    is isomorphic to the truncation triangle with respect to $\psi$. The hypothesis on the Thomason filtration tells us that $\psi(i)$ is contained in $\phi(i)$ for all $i$. This ensures that $\mathcal{U}_{\psi}$ is contained in $\mathcal{U}_\phi$, and so, $\mathcal{V}_\phi[-1]$ is contained in $\mathcal{V}_\psi [-1]$. Hence, $\tau^{\geq 1}_\phi F$ belongs to $\mathcal{V}_\psi [-1]$. It follows from \Cref{lemma:AJS2010_4.2_generalization } that $\tau_\phi^{\leq 0} D^{\geq j}_{\operatorname{qc}}(X)$ and $\tau_\phi^{\geq 1} D^{\geq j}_{\operatorname{qc}}(X)$ are both contained in $D^{\geq j}_{\operatorname{qc}}(X)$ for all integers $j$. Then $\tau^{\leq 0}_\phi F$ belongs to $\mathcal{U}_\psi$, and so, $\tau^{\leq 0}_\phi  F$ is isomorphic to $\tau^{\leq 0}_\psi F$ as desired.
\end{proof}

\begin{corollary}\label{cor:filtrations_agreeing_after_some_point_have_same_truncation}
    Let $i,j$ be integers. Suppose $\phi,\psi$ are Thomason filtrations on $X$ such that $\phi (i) = \psi(i)$ for all $i\geq j$. For $F\in D^{\geq j}_{\operatorname{qc}}(X)$, one has $\tau^{\leq 0}_\phi F\cong \tau^{\leq 0}_\psi F,\tau^{\geq 1}_\phi F\cong \tau^{\geq 1}_\psi F \in D_{\operatorname{qc}}(X)$.
\end{corollary}

\begin{proof}
    This follows from \Cref{prop:filtrations_agreeing_after_some_point} as $\phi^{\geq j} = \psi^{\geq j}$.
\end{proof}

\begin{corollary}\label{cor:cohomology_for_filtrations_agreeing_on_strip}
    Assume $X=\operatorname{Spec}(R)$. Let $a\leq b$ be integers. Suppose $\phi$ and $\psi$ are two Thomason filtrations on $X$ where $\phi (i) = \psi(i)$ for all $i$ in $[a,b]\cap \mathbb{Z}$. Then, for any object $F$ in $D^{\geq a}_{\operatorname{qc}}(R)$, and integer $i$ in $[a,b]\cap \mathbb{Z}$, one has that $H^i(\tau^{\leq 0}_\phi F)\cong H^i(\tau^{\leq 0}_\psi F)$.
\end{corollary}

\begin{proof}
    The hypothesis ensures that $\phi \cap ( \phi _{\operatorname{st}}[-b])(t)= \psi \cap( \phi _{\operatorname{st}}[-b])(t)$ if $t\geq a$. It follows from \cite[Lemma 5.2]{Smith:2022} and \Cref{cor:filtrations_agreeing_after_some_point_have_same_truncation} that $\tau^{\leq b}_{\operatorname{st}}\circ \tau^{\leq 0}_\phi F$ is isomorphic to $\tau^{\leq b}_{\operatorname{st}}\circ \tau^{\leq 0}_{\psi} F$. However, $\tau^{\leq b}_{\operatorname{st}}$ does not alter cohomology in degrees less than or equal to $b$, so the claim follows.
\end{proof}

\begin{lemma}\label{lem:all_associated_t_structures_restrict}
    Let $\phi$ be a Thomason filtration on $X$ and $Z$ be a closed subset of $X$. If $E$ is an object of $D_{\operatorname{qc},Z}(X)$, then both $\tau^{\leq 0}_\phi E$ and $\tau^{\geq 1}_\phi E$ belong to $D_{\operatorname{qc},Z}(X)$. Consequently, the $t$-structure associated to $\phi$ on $D_{\operatorname{qc}}(X)$ restricts to $D_{\operatorname{qc},Z}(X)$ and is tensor on $D_{\operatorname{qc},Z}(X)$.
\end{lemma}

\begin{proof}
    This follows from the fact that passing to stalks commutes with the connective and coconnective truncations associated to $\phi$. Indeed, it can be observed by \cite[Proposition 1.13]{AlonsoTarrio/JeremiasLopez/Saorin:2010} and \cite[Proposition 3.5]{Lank:2025}. That the last claim holds follows from the fact $D^{\leq 0}_{\operatorname{qc},Z}(X)\subseteq D^{\leq 0}_{\operatorname{qc}}(X)$.
\end{proof} 

\begin{proposition}\label{prop:truncations_preserve_support}
    Let $Z$ be a closed subset of $X$. Consider a Thomason filtration $\phi$ on $X$ and an object $F$ of $D^{\geq j}_{\operatorname{qc},Z}(X)$. For $\psi = \phi\cap Z$, one has $\tau^{\leq 0}_\phi F\cong \tau^{\leq 0}_\psi F$ and $\tau^{\geq 1}_\phi F\cong \tau^{\geq 1}_\psi F$ in $D_{\operatorname{qc}}(X)$.
\end{proposition}

\begin{proof}
    It suffices to show that $\tau^{\leq 0}_\phi F$ is an object of $\mathcal{U}_\psi$. Note that $\operatorname{supp}(H^i(\tau^{\leq 0}_{\phi}(A)))$ is contained in $\phi(i)$ for all integers $i$. Recall that for any object $A$ in $D_{\operatorname{qc}}(X)$ we have that $\operatorname{supp}(H^i(A))$ is contained in $\operatorname{supp}(A)$. Then $\operatorname{supp}(H^i(\tau^{\leq 0}_{\phi}(F))) \subseteq \phi (i)\cap Z$ as desired. The desired claim follows from \cite[Theorem 3.6]{Lank:2025} and \cite[Example 2.1]{Hrbek/Lank/LeGros/Pavon:2026}.
\end{proof}

\begin{corollary}\label{cor:automatic_restriction_to_smaller_Thomason}
    Let $W,Z$ be closed subsets of $X$ where $W\subseteq Z$. If $\phi$ is a Thomason filtration on $X$ whose associated $t$-structure on $D_{\operatorname{qc}}(X)$ restricts to $\operatorname{Perf}_Z(X)$, then the associated $t$-structure of $\phi$ also restricts to $\operatorname{Perf}_W(X)$. In particular, the restricted $t$-structure is tensor on $\operatorname{Perf}_W (X)$.
\end{corollary}

\begin{proof}
    That the first claim follows from \Cref{lem:all_associated_t_structures_restrict} and \Cref{prop:truncations_preserve_support}. The last claim holds by the fact that $D^{\leq 0}_{\operatorname{qc},Z}(X)\subseteq D^{\leq 0}_{\operatorname{qc}}(X)$.
\end{proof}

\begin{lemma}\label{lem:thomason_closed_aisle_lift}
    Let $Z$ be a closed subset of $X$. Suppose $\mathcal{U}$ is an aisle on $D_{\operatorname{qc},Z}(X)$ such that $\mathcal{U} = \overline{ \langle \mathcal{S} \rangle}^{(-\infty,0]}$ where $\mathcal{S}\subseteq \operatorname{Perf}_Z (X)$. Then $\mathcal{U}^\prime:= \overline{ \langle \mathcal{S} \rangle}^{(-\infty,0]}_{\otimes}$ on $D_{\operatorname{qc}}(X)$ is associated to a Thomason filtration $\phi$ on $X$ such that $\phi (n)\subseteq Z$ for all $n$.
\end{lemma}

\begin{proof}
    By \cite[Theorem 3.6]{Lank:2025}, the associated Thomason filtration of $\mathcal{U}^\prime$ is given by the following:
    \begin{displaymath}
        \phi_{\mathcal{U}^\prime} (i)= \bigcup_{P \in \mathcal{S}} \bigcup_{j\geq i} \operatorname{supp}(\mathcal{H}^j (P)).
    \end{displaymath}
    However, every object $P$ in $\mathcal{S}$ is supported in $Z$, and so the desired claim follows.
\end{proof}

\begin{proposition}\label{prop:tensor_t_structures_on_D_Z}
    Let $Z$ be a closed subset of $X$. Then there is a bijective correspondence between the collection of compactly generated tensor $t$-structures on $D_{\operatorname{qc},Z}(X)$ and the collection of Thomason filtrations $\phi$ on $X$  such that $\phi (n)$ contained in $Z$ for all $n$. 
\end{proposition}

\begin{proof}
    This follows from \Cref{lem:all_associated_t_structures_restrict} and \Cref{lem:thomason_closed_aisle_lift} if coupled with \cite[Theorem 3.6]{Lank:2025} and \cite[Example 2.1]{Hrbek/Lank/LeGros/Pavon:2026}. Specifically, a compactly generated $\otimes$-aisle $\mathcal{A}$ on $D_{\operatorname{qc},Z}(X)$ is assigned the Thomason filtration corresponding to the aisle $\overline{ \langle \mathcal{A} \rangle}^{\leq}_{\otimes}$ on $D_{\operatorname{qc}}(X)$. And a Thomason filtration $\phi$ on $X$ supported on $Z$ is assigned the compactly generated $\otimes$-aisle $\mathcal{U}_\phi \cap D_{\operatorname{qc}}(X)$. This gives us maps for the correspondence. A direct argument will verify this is a one-to-one correspondence.
\end{proof}

\begin{corollary}\label{cor:lifitng_tensor_structure_from_perfect_categories}
    Let $Z$ be a closed subset of $X$. Any tensor $t$-structure on $\operatorname{Perf}_Z (X)$ is the restriction of a $t$-structure on $D_{\operatorname{qc},Z}(X)$ that is associated to Thomason filtration $\phi$ on $X$ such that $\phi (n)$ is contained in $Z$ for all $n$.
\end{corollary}

\begin{proof}
    First we show that any tensor $t$-structure on $\operatorname{Perf}_Z (X)$ is the restriction of a compactly generated tensor $t$-structure on $D_{\operatorname{qc},Z}(X)$. Let $\mathcal{A}$ be a $\otimes$-aisle on $\operatorname{Perf}_Z (X)$. We claim that $\operatorname{Coprod}(\mathcal{A}) = \overline{ \langle \mathcal{A} \rangle}^{(-\infty,0]}_{\otimes}$. It suffices to show $D^{\leq 0}_{\operatorname{qc}}(X) \otimes^{\mathbf{L}} \mathcal{A} \subseteq \mathcal{A}$. By \Cref{lem:lifitng_tensor_structure_from_perfect_categorie_abstract}, $\operatorname{Coprod}(\mathcal{A}) = \overline{ \langle \mathcal{A} \rangle}^{(-\infty,0]}$ is a compactly generated aisle on $D_{\operatorname{qc},Z}(X)$ that restricts to $\mathcal{A}$. Observe that $\operatorname{Coprod}(\mathcal{A})$ is also an aisle on $D_{\operatorname{qc}}(X)$. Indeed, $\operatorname{Coprod}(-)$ computed for $\mathcal{A}$ in $D_{\operatorname{qc},Z}(X)$ or $D_{\operatorname{qc}}(X)$ gives the same category as $D_{\operatorname{qc},Z}(X)$ is a strictly full subcategory of $D_{\operatorname{qc}}(X)$.
    
    It follows from \cite[Theorem 3.6]{Lank:2025}, \cite[Example 2.1]{Hrbek/Lank/LeGros/Pavon:2026}, and \Cref{prop:tensor_t_structures_on_D_Z} that
    \begin{displaymath}
        D^{\leq 0}_{\operatorname{qc},Z}(X) = \operatorname{Coprod}(\operatorname{Perf}_Z (X) \cap D^{\leq 0}_{\operatorname{qc},Z}(X)).
    \end{displaymath}
    From $\mathcal{A}$ being an $\otimes$-aisle on $\operatorname{Perf}_Z (X)$, it follows that $(\operatorname{Perf}_Z (X)\cap D^{\leq 0}_{\operatorname{qc}}(X))\otimes^{\mathbf{L}}\mathcal{A}\subseteq \mathcal{A}$. Let $E\in D^{\leq 0}_{\operatorname{qc},Z}(X)$ and $A\in \mathcal{A}$. By \cite[Theorem 12.1]{Keller/Nicolas:2013}, there exists a distinguished triangle in $D_{\operatorname{qc},Z}(X)$ of the form
    \begin{displaymath}
        \coprod_{i\geq 0} E_i \to E \to \coprod_{i\geq 0} E_i[1] \to (\coprod_{i\geq 0} E_i)[1]
    \end{displaymath}
    where each $E_i$ is an $i$-fold extension of small coproducts of non-negative shifts of objects in $\operatorname{Perf}_Z (X)\cap D^{\leq 0}_{\operatorname{qc}}(X)$. Recall that $-\otimes^{\mathbf{L}} A$ is left adjoint to $\operatorname{\mathbf{R}\mathcal{H}\! \mathit{om}}(A,-)$. So we have an isomorphism $(\coprod_{i\geq 0} E_i)\otimes^{\mathbf{L}} A \cong \coprod_{i\geq 0} (E_i\otimes^{\mathbf{L}} A)$. Our hypothesis ensures each $E_i \otimes^{\mathbf{L}} A\in \mathcal{A}$. Hence, we see that $\coprod_{i\geq 0} E_i\otimes^{\mathbf{L}} A\in \operatorname{Coprod}(\mathcal{A})$, and so $E\otimes^{\mathbf{L}} A\in \operatorname{Coprod}(\mathcal{A})$ being that aisles are closed under extensions. This shows that $D^{\leq 0}_{\operatorname{qc}}(X) \otimes^{\mathbf{L}} A\in \operatorname{Coprod}(\mathcal{A})$. 

    Now we prove the desired claim. It was shown above that $\operatorname{Coprod}(\mathcal{A})$ is a compactly generated $\otimes$-aisle on $D_{\operatorname{qc}}(X)$. Then \Cref{prop:tensor_t_structures_on_D_Z} gives us the desired associated Thomason filtration $\phi$ on $X$.
\end{proof}

\begin{lemma}\label{lem:induced_t_structure_verdier_localization1}
    Suppose $j^\ast \colon \mathcal{T} \to \mathcal{S}$ is an essentially surjective exact functor between triangulated categories. For any $t$-structure $(\mathcal{A}, \mathcal{B})$ on $\mathcal{T}$, the following are equivalent:
    \begin{enumerate}
        \item $(j^\ast \mathcal{A}, j^\ast \mathcal{B})$ is $t$-structure on $\mathcal{S}$
        \item $\operatorname{Hom}(j^\ast \mathcal{A},j^\ast \mathcal{B}[-1])=0$.
    \end{enumerate}
\end{lemma}

\begin{proof}
    The forward direction follows because the orthogonality conditions in the definition of a $t$-structure, whereas the converse direction follows from the fact that $j^\ast$ is an essentially surjective exact functor. 
\end{proof}

\begin{lemma}\label{lem:induced_t_structure_verdier_localization2}
    With the notation of \Cref{lem:induced_t_structure_verdier_localization1} where each triangulated category is compactly generated. Suppose $j^\ast$ admits a fully faithful exact right adjoint $j_\ast \colon \mathcal{S} \to \mathcal{T}$ which preserves small coproducts and $(j^\ast \mathcal{A}, j^\ast \mathcal{B})$ is a $t$-structure on $\mathcal{S}$. If $(\mathcal{A},\mathcal{B})$ is compactly generated, then so is $(j^\ast \mathcal{A}, j^\ast \mathcal{B})$.
\end{lemma}

\begin{proof}
    Let $\mathcal{P}$ be any collection of compact objects of $\mathcal{T}$ that generate $\mathcal{A}$. By \cite[Theorem 5.1]{Neeman:1996}, $j_\ast$ preserving small coproducts implies $j^\ast$ preserves compacts. Hence, $j^\ast \mathcal{P}$ is a collection of compact objects contained in the aisle $j^\ast \mathcal{A}$. Note that every object in $j^\ast \mathcal{A}$ is isomorphic to an object of the form $j^\ast A$ for $A$ in $\mathcal{A}$, and similarly for $j^\ast \mathcal{P}$. Suppose $j^\ast E$ is an object of $j^\ast \mathcal{A}$ such that $\operatorname{Hom}(j^\ast P,j^\ast E)=0$ for all $P$ in $\mathcal{P}$. It follows from adjunction that $\operatorname{Hom}(P, j_\ast j^\ast E)=0$ for all $P$ in $\mathcal{P}$. Let $A$ be an object of $\mathcal{A}$. From \cite[Theorem 12.1]{Keller/Nicolas:2013}, there exists a distinguished triangle in $\mathcal{T}$
    \begin{displaymath}
        \coprod_{i\geq 0} A_i \to A \to \coprod_{i\geq 0} A_i[1] \to (\coprod_{i\geq 0} A_i)[1]
    \end{displaymath}
    where each $A_i$ is an $i$-fold extension of small coproducts of non-negative shifts of objects in $\mathcal{P}$. It follows from applying the cohomological functor $\operatorname{Hom}(-,j_\ast j^\ast E) \colon \mathcal{T} \to \operatorname{Ab}$ to this distinguished triangle that $\operatorname{Hom}(A,j_\ast j^\ast E)=0$. This tells us that $j^\ast E$ belongs to $j^\ast \mathcal{B}[-1]$ because $j^\ast j_\ast j^\ast E$ is isomorphic to $j^\ast E$ and $j_\ast j^\ast E$ is in $\mathcal{B}[-1]$. Then $j^\ast E [1]$ is an object of $j^\ast \mathcal{B}$. However, $j^\ast E[1]$ is in $j^\ast \mathcal{A}$ and $\operatorname{Hom}(j^\ast \mathcal{A},j^\ast \mathcal{B})=0$, which tells us $j^\ast E$ must be the zero object. Therefore, $j^\ast \mathcal{P}$ compactly generates the aisle $j^\ast \mathcal{A}$. Indeed, it is clear that $j^\ast \mathcal{P}\subseteq j^\ast \mathcal{A}$. By \cite[Theorem 2.3]{Neeman:2021a}, $j^\ast \mathcal{P}$ is an aisle. Let $E\in j^\ast \mathcal{A}$. Consider the truncation triangle of $E$ with respect to $j^\ast \mathcal{P}$,
    \begin{displaymath}
        \tau^{\leq 0} E \to E \to \tau^{\geq 1} E \to (\tau^{\leq 0} E)[1].
    \end{displaymath}
    Since $j^\ast \mathcal{A}$ is an aisle, $\tau^{\geq 1} E \in j^\ast \mathcal{A}$. To see, use that $\tau^{\geq 1} E$ is an extension of $E,(\tau^{\leq 0} E)[1]\in j^\ast \mathcal{A}$. However, $\tau^{\geq 1} E \in (j^\ast \mathcal{P})^\perp$, so $\tau^{\geq 1} E \cong 0$, which completes the proof.
\end{proof}

\begin{lemma}\label{lem:stalks_to_global_tensor_compatible}
    Let $(\mathcal{A},\mathcal{B})$ be a $t$-structure on $D_{\operatorname{qc}}(X)$. If $(\mathcal{A}_p,\mathcal{B}_p)$ is a $t$-structure on $D_{\operatorname{qc}}(\mathcal{O}_{X,p})$ for all $p$ in $X$, then $(\mathcal{A},\mathcal{B})$ is a tensor $t$-structure on $D_{\operatorname{qc}}(X)$.
\end{lemma}

\begin{proof}
    Recall that all $t$-structures on $D_{\operatorname{qc}}(W)$ are tensor compatible for $W$ an affine Noetherian scheme (see \cite[Corollary 3.8]{Lank:2025}). We freely use this fact in the proof. Let $A$ be an object in $\mathcal{A}$. Choose an object $E$ in $D^{\leq 0}_{\operatorname{qc}}(X)$. There exists a distinguished triangle in $D_{\operatorname{qc}}(X)$:
    \begin{displaymath}
        A^\prime \xrightarrow{f} E\overset{\mathbf{L}}{\otimes} A \to B \to A^\prime[1]
    \end{displaymath}
    where $A^\prime$ is in $\mathcal{A}$ and $B$ is in $\mathcal{B}[-1]$. We know that $(\mathcal{A}_p, \mathcal{B}_p)$ is tensor $t$-structure on $D_{\operatorname{qc}}(\mathcal{O}_{X,p})$ for all points $p$ in $X$. Hence, there exists a distinguished triangle in $D_{\operatorname{qc}}(X)$:
    \begin{displaymath}
        A^\prime_p \xrightarrow{f} (E\overset{\mathbf{L}}{\otimes} A)_p \to B_p \to A^\prime_p[1]
    \end{displaymath}
    where $A^\prime_p$ is in $\mathcal{A}_p$ and $B_p$ is in $\mathcal{B}_p[-1]$. It follows that $E_p \overset{\mathbf{L}}{\otimes} A_p$ is in $\mathcal{A}_p$ because $E_p$ is an object of $D_{\operatorname{qc}}^{\leq 0}(\mathcal{O}_{X,p})$. However, $E_p \overset{\mathbf{L}}{\otimes} A_p$ is isomorphic to $(E \overset{\mathbf{L}}{\otimes} A)_p$ in $D_{\operatorname{qc}}(\mathcal{O}_{X,p})$. Indeed, it follows from the fact that localization is flat. This observation coupled with the distinguished triangle above implies $B_p$ vanishes at all $p$ in $X$. Hence, $B$ has empty support, so it must be the zero object, ensuring that  $A^\prime$ is isomorphic to $E \overset{\mathbf{L}}{\otimes} A$ in $D_{\operatorname{qc}}(X)$. Thus, $E\overset{\mathbf{L}}{\otimes} A$ belongs to $\mathcal{A}$, telling us this is a $\otimes$-aisle, which completes the proof. 
\end{proof}

\begin{proposition}\label{prop:tensor_compatible_local_to_global}
    Let $Z$ be a closed subset of $X$. Suppose $(\mathcal{A},\mathcal{B})$ is a compactly generated $t$-structure on $D_{\operatorname{qc}}(X)$. Then the following are equivalent:
    \begin{enumerate}
        \item $(\mathcal{A},\mathcal{B})$ restricts to a tensor $t$-structure on $D_{\operatorname{qc},Z}(X)$
        \item $(j^\ast \mathcal{A}, j^\ast \mathcal{B})$ restricts to a compactly generated $t$-structure on $D_{\operatorname{qc},Z\cap U}(U)$ for all open immersions $j \colon U \to X$ where $U$ is affine
        \item $(\mathcal{A}_p,\mathcal{B}_p)$ restricts to a compactly generated $t$-structure on $D_{\operatorname{qc},Z\cap \operatorname{Spec}(\mathcal{O}_{X,p})}(\mathcal{O}_{X,p})$ for all $p$ in $X$.
    \end{enumerate}
\end{proposition}

\begin{proof}
    This is known to experts but we include it for convenience. It follows from \Cref{lem:all_associated_t_structures_restrict} and \Cref{prop:tensor_t_structures_on_D_Z} that we only need to check the case where $Z=X$. On one hand, $(3)\implies (1)$ is \Cref{lem:stalks_to_global_tensor_compatible}\footnote{Note that being compactly generated is not required for this direction}, whereas one the other hand,  $(2)\implies (3)$ follows from \cite[Proposition 2.9 \& Theorem 3.10]{AlonsoTarrio/JeremiasLopez/Saorin:2010}. We verify $(1)\implies (2)$. By \cite[Proposition 3.5]{Lank:2025}, if $j \colon U \to X$ is an open immersion where $U$ is affine, then $\operatorname{Hom}(j^\ast A, j^\ast B)=0$ for any pair of objects $A$ in $\mathcal{A}$ and $B$ in $\mathcal{B}$. It follows from \Cref{lem:induced_t_structure_verdier_localization2} that $(j^\ast \mathcal{A}, j^\ast \mathcal{B})$ is a compactly generated $t$-structure on $D_{\operatorname{qc}}(U)$.
\end{proof}

\begin{lemma}\label{lem:relative_d^b_coh_verdier_localization}
    Let $X$ be a Noetherian scheme. Suppose $Z$ is a closed subset of $X$. Then there are Verdier localizations:
    \begin{enumerate}
        \item $j^\ast \colon D^b_{\operatorname{coh},Z}(X) \to D^b_{\operatorname{coh},Z\cap U}(U)$ where $j\colon U \to X$ is an open immersion from an affine scheme
        \item $s^\ast \colon D^b_{\operatorname{coh},Z}(X) \to D^b_{\operatorname{coh},Z\cap \operatorname{Spec}(\mathcal{O}_{X,p})}(\mathcal{O}_{X,p})$ where $s\colon \operatorname{Spec}(\mathcal{O}_{X,p}) \to X$ is the natural morphism for $p$ in $X$.
    \end{enumerate}
\end{lemma}

\begin{proof}[Proof sketch]
    It follows from \cite[Proposition 3.9.2]{Lipman/Hashimoto:2009} that objects with bounded quasi-coherent cohomology have image under both derived pushforwards also with bounded quasi-coherent cohomology\footnote{This can also be inferred, essentially in the same manner, from \cite[\href{https://stacks.math.columbia.edu/tag/08D5}{Tag 08D5}]{StacksProject} and \cite[\href{https://stacks.math.columbia.edu/tag/05T6}{Tag 05T6}]{StacksProject}.}. Moreover, both functors send objects to those in $D_{\operatorname{qc}}(X)$ which are supported in $Z$. It is evident that both associated derived pullback functors are well-defined as stated in the claim. The desired claims follow by arguing in a similar fashion to \cite[Theorem 4.4]{Elagin/Lunts/Schnurer:2020}.
\end{proof}

\begin{corollary}\label{cor:descent_for_thomason_on_bounded}
    Let $X$ be a Noetherian scheme. Let $Z$ be a closed subset of $X$. Suppose $\phi$ is a Thomason filtration on $X$ whose associated $t$-structure $(\mathcal{A},\mathcal{B})$ on $D_{\operatorname{qc}}(X)$ restricts to a $t$-structure on $D^b_{\operatorname{coh},Z}(X)$. 
    \begin{enumerate}
        \item If $j \colon U \to X$ is an open immersion from an affine scheme, then the $t$-structure associated to the Thomason filtration $\phi$ on $D_{\operatorname{qc},Z\cap U}(U)$ restricts to $D^b_{\operatorname{coh},Z\cap U}(U)$.
        \item For any point $p$ in $X$, the $t$-structure associated to the Thomason filtration $\phi_p$ on $D_{\operatorname{qc}}(\mathcal{O}_{X,p})$ restricts to $D^b_{\operatorname{coh},Z_p}(\mathcal{O}_{X,p})$ where $Z_p:=Z\cap \operatorname{Spec}(\mathcal{O}_{X,p})$.
    \end{enumerate}
\end{corollary}

\begin{proof}
    It suffices by \Cref{lem:relative_d^b_coh_verdier_localization} that we only need to check $(1)$ as $(2)$ can be shown in a similar manner. Let $j \colon U \to X$ be an open immersion from an affine scheme. We know from \Cref{prop:tensor_compatible_local_to_global} that $(j^\ast \mathcal{A}, j^\ast \mathcal{B})$ is a compactly generated $t$-structure on $D_{\operatorname{qc},Z\cap U}(U)$. Fix an object $P$ in $D^b_{\operatorname{coh},Z\cap U}(U)$. There exists $P^\prime$ in $D^b_{\operatorname{coh},Z}(X)$ such that $P$ is isomorphic to $j^\ast P^\prime$ in $D^b_{\operatorname{coh},Z\cap U} (U)$ by \Cref{lem:relative_d^b_coh_verdier_localization}. We can find a distinguished triangle in $D_{\operatorname{qc},Z}(X)$:
    \begin{displaymath}
        A \to P^\prime \to B \to A[1]
    \end{displaymath}
    where $A$ is in $\mathcal{A}$ and $B$ is in $\mathcal{B}[-1]$. However, this distinguished triangle must be the corresponding truncation triangle of $P^\prime$ for the $t$-structure $(\mathcal{A}\cap D^b_{\operatorname{coh},Z} (X), \mathcal{B}\cap D^b_{\operatorname{coh},Z} (X))$ on $D^b_{\operatorname{coh},Z}(X)$, and so $A,B$ belong to $D^b_{\operatorname{coh},Z}(X)$. The restriction of $\phi$ to $U$ has the associated $t$-structure $(j^\ast \mathcal{A},j^\ast \mathcal{B})$, and so, the corresponding truncation triangle for $j^\ast P$ is the following:
    \begin{displaymath}
        j^\ast A \to j^\ast P^\prime \to j^\ast B \to j^\ast A[1]
    \end{displaymath}
    where $j^\ast A$ is in $j^\ast \mathcal{A}$ and $j^\ast B$ is in $j^\ast \mathcal{B}[-1]$. However, we know that $j^\ast A, j^\ast B$ are in $D^b_{\operatorname{coh}, Z\cap U}(U)$, ensuring that $j^\ast P^\prime$ belongs to $D^b_{\operatorname{coh}, Z\cap U}(U)$. This completes the proof.
\end{proof}

\begin{lemma}\label{lem:local_to_global_aisle}
    Let $Z$ be a closed subset of $X$. Suppose $(\mathcal{A},\mathcal{B})$ is a $t$-structure on $D_{\operatorname{qc},Z}(X)$ such that for all open immersions $j \colon U \to X$ where $U$ is affine, the pair $(j^\ast \mathcal{A}, j^\ast \mathcal{B})$ is $t$-structure on $D_{\operatorname{qc},Z\cap U}(U)$. An object $E$ in $D_{\operatorname{qc},Z}(X)$ belongs to $\mathcal{A}$ if, and only if, $j^\ast E$ is in $j^\ast \mathcal{A}$ for all open immersions $j \colon U \to X$ where $U$ is affine.
\end{lemma}

\begin{proof}
    The forward direction holds from our hypothesis, so we check the converse. Let $E$ be an object of $D_{\operatorname{qc},Z}(X)$ such that $j^\ast E$ is in $j^\ast \mathcal{A}$ for all open immersions $j \colon U \to X$ where $U$ is affine. There exists a distinguished triangle in $D_{\operatorname{qc}}(X)$:
    \begin{displaymath}
        A \xrightarrow{f} E \to B \to A[1]
    \end{displaymath}
    where $A$ is in $\mathcal{A}$ and $B$ is in $\mathcal{B}[-1]$. From our hypothesis, we have the corresponding distinguished triangle of $j^\ast E$ for the $t$-structure $(j^\ast \mathcal{A},j^\ast \mathcal{B})$ in $D_{\operatorname{qc}, Z\cap U}(U)$ is given by:
    \begin{displaymath}
        j^\ast A \xrightarrow{j^\ast ( f )} j^\ast E \to j^\ast B \to j^\ast A[1].
    \end{displaymath}
    The hypothesis that $j^\ast E$ belongs to $j^\ast \mathcal{A}$ tells us that the map $j^\ast (f) \colon j^\ast A \to j^\ast E$ is a isomorphism in $D_{\operatorname{qc}, Z\cap U}(U)$, which tells us $j^\ast B$ is zero. However, it follows that $j^\ast B$ is zero for all open immersions $j\colon U \to X$ with $U$ affine, and so $B$ must be zero in $D_{\operatorname{qc},Z}(X)$. Hence, $f\colon A \to E$ is a isomorphism, ensuring that $E$ belongs to $\mathcal{A}$ as desired.
\end{proof}

\begin{lemma}\label{lem:global_to_open_perfect_supported}
    Let $Z$ be a closed subset of $X$. Then there are essentially dense functors:
    \begin{enumerate}
        \item $j^\ast \colon \operatorname{Perf}_Z (X)\to \operatorname{Perf}_{Z\cap U} (U)$ where $j\colon U \to X$ is an open immersion
        \item $s^\ast \colon \operatorname{Perf}_Z (X) \to \operatorname{Perf}_{Z\cap \operatorname{Spec}(\mathcal{O}_{X,p})} (\mathcal{O}_{X,p})$ where $s\colon \operatorname{Spec}(\mathcal{O}_{X,p}) \to X$ is the natural morphism for $p$ in $X$.
    \end{enumerate}
\end{lemma}

\begin{proof}
    There are Verdier localizations induced by both morphisms:
    \begin{displaymath}
        \begin{aligned}
            & D_{\operatorname{qc},Z\setminus (U \cap Z)}(X)  \to D_{\operatorname{qc},Z} (X) \xrightarrow{j^\ast} D_{\operatorname{qc},Z \cap U} (U), 
            \\& D_{\operatorname{qc},Z\setminus (\{ q\in X : q\not\in \overline{\{ p \}}\} \cap Z)}(X) \to D_{\operatorname{qc},Z} (X) \xrightarrow{s^\ast} D_{\operatorname{qc},Z \cap \operatorname{Spec}(\mathcal{O}_{X,p})} (\mathcal{O}_{X,p}).
        \end{aligned}
    \end{displaymath}
    Then the desired claim follows by \cite[Theorem 2.1.4]{Neeman:1996}.
\end{proof}

\begin{proposition}\label{prop:reduction_to_affines_stalks}
   Let $Z$ be a closed subset of $X$. Suppose $\phi$ is a Thomason filtration on $X$ supported in $Z$. Assume the $t$-structure on $D_{\operatorname{qc}}(X)$ associated to $\phi$ restricts to a $t$-structure on $\operatorname{Perf}_Z (X)$. 
   \begin{enumerate}
       \item If $j \colon U \to X$ is an open immersion from an affine scheme, then the $t$-structure associated to the Thomason filtration $\phi\cap U$ on $D_{\operatorname{qc}}(U)$ restricts to $\operatorname{Perf}_{Z\cap U} (U)$.
       \item For any point $p$ in $Z$, the $t$-structure associated to the Thomason filtration $\phi_p$ on $D_{\operatorname{qc}}(\mathcal{O}_{X,p})$ restricts to $\operatorname{Perf}_{Z\cap \operatorname{Spec}(\mathcal{O}_{X,p})}(\mathcal{O}_{X,p})$.
   \end{enumerate}
\end{proposition}

\begin{proof}
   We will prove $(1)$ as $(2)$ is shown in a similar manner. Choose an object $P$ in $\operatorname{Perf}_{Z\cap U} (U)$. From \Cref{lem:global_to_open_perfect_supported}, there exists $P^\prime$ in $\operatorname{Perf}_Z (X)$ such that $P$ is a direct summand of $j^\ast P^\prime$ in $\operatorname{Perf}_{Z\cap U}( U)$. Consider the truncation triangle in $D_{\operatorname{qc}}(X)$:
   \begin{displaymath}
       A \to P^\prime \to B \to A[1]
   \end{displaymath}
   where $A$ is in $\mathcal{A}$ and $B$ is in $\mathcal{B}[-1]$. Our hypothesis tells us that $A,B$ belong to $\operatorname{Perf}_Z (X)$ as $P^\prime$ is in $\operatorname{Perf}_Z (X)$ and $\mathcal{A}$ restricts to $\operatorname{Perf}_Z(X)$. The restriction of $\phi$ to $D_{\operatorname{qc}}(U)$ has associated $t$-structure $(j^\ast \mathcal{A},j^\ast \mathcal{B})$, and we have a corresponding truncation triangle for $j^\ast P^\prime$:
   \begin{displaymath}
       j^\ast A \to j^\ast P^\prime \to j^\ast B \to j^\ast A[1].
   \end{displaymath}
   where $j^\ast A$ is in $j^\ast \mathcal{A}$ and $j^\ast B$ is in $j^\ast \mathcal{B}[-1]$, both belonging to $\operatorname{Perf}_{Z\cap U} (U)$. There exist maps $f\colon P \to j^\ast P^\prime$ and $g\colon j^\ast P^\prime \to P$ whose composition is the identity. Denote the truncation functors of the $t$-structure $(j^\ast \mathcal{A},j^\ast \mathcal{B})$ on $D_{\operatorname{qc}}(U)$ by $\tau^{\leq 0} \colon D_{\operatorname{qc}}(U) \to \mathcal{A}$ and $\tau^{\geq 1} \colon D_{\operatorname{qc}}(U) \to \mathcal{B}$. Applying the connective truncation functors $\tau^{\leq 0}, \tau^{\geq 1}$ tells us that $\tau^{\leq 0} P$ is a direct summand of $j^\ast A$ and $\tau^{\geq 1} P$ is a direct summand of $j^\ast B$. In other words, the truncations of $P$ with respect to $(j^\ast \mathcal{A}, j^\ast \mathcal{B})$ belong to $\operatorname{Perf}_{Z\cap U} (U)$.
\end{proof}

\section{Relative weak Cousin}
\label{sec:results1}

Let $X$ be a Noetherian scheme. We need a relative version of a notion in \cite[\S 4]{AlonsoTarrio/JeremiasLopez/Saorin:2010}.

\begin{definition}\label{def:weak_cousin_relative}
    Let $Z$ be a subset of $X$. A Thomason filtration $\phi$ on $X$ is said to be \textbf{weak Cousin across $Z$} if for any direct generalization  $p \rightsquigarrow q$ with $p,q$ in $Z$, if $q$ is in $\phi (n)$, then $p$ is in $\phi (n-1)$.
\end{definition}

\begin{remark}\label{rmk:weak_cousin_across_to_literal_weak_cousin_on_intersection}
    Let $Z$ be a subset of $X$. A Thomason filtration $\phi$ on $X$ is weak Cousin across $Z$ if, and only if, $\phi\cap Z$ is weak Cousin on $Z$.
\end{remark}

\begin{lemma}\label{lem:weak_cousin_iff_local}
    Let $Z$ be a closed subset of $X$. If $\phi \colon \mathbb{Z} \to \operatorname{Spcl}(X)$ is Thomason filtration on $X$, then the following are equivalent:
    \begin{enumerate}
        \item $\phi$ satisfies the weak Cousin condition across $Z$
        \item $\phi\cap U$ satisfies the weak Cousin condition across $Z \cap U$ for any affine open subscheme $U$ of $X$ where $Z\cap U\not= \emptyset$
        \item $\phi_s$ satisfies the weak Cousin condition across $Z\cap \operatorname{Spec}(\mathcal{O}_{X,s})$ for all $s$ in $Z$.
    \end{enumerate}
\end{lemma}

\begin{proof}
    To start, we show $(1)\implies (2)$. Suppose $\phi$ satisfies the weak Cousin condition across $Z$. Let $j \colon U \to X$ be an open immersion such that $Z\cap U\not= \emptyset$. Fix an integer $n$. Consider a direct generalization  $p \rightsquigarrow q$ with $q$ in $( \phi \cap U)(n)$ and $p$ in $Z\cap U$. Then $q$ is in $\overline{\{p  \}}$ (i.e.\ the closure of $p$) in $X$, which ensures that $p$ is in $\phi (n-1)$, and so, $p$ belongs to $( \phi \cap U)(n-1)$. 

    Next, we show $(2)\implies (3)$. Let $s$ be a point of $Z$. Choose an affine open neighbourhood $U$ of $s$ in $X$. Fix an integer $n$. Suppose $p \rightsquigarrow q$ is a direct generalization  in $\operatorname{Spec}(\mathcal{O}_{X,s})$ with $p,q$ in $Z \cap \operatorname{Spec}(\mathcal{O}_{X,s})$ and $q$ is in $\phi_s(n)$. It follows that this gives us a direct generalization  of points in $Z\cap U$ as open sets are generalisation closed and $s$ is in $U$. Thus, the hypothesis implies $p$ is in $( \phi \cap U)(n-1)$, ensuring that $p$ belongs to $\phi_s(n-1)$ as desired.

    Lastly, we show $(3)\implies (1)$. Fix an integer $n$. Let $p \rightsquigarrow q$ be a direct generalization  in $Z$ where $q$ is in $\phi(n)$. Then $p \rightsquigarrow q$ is a direct generalization  in $\operatorname{Spec}(\mathcal{O}_{X,q})$ where $q$ is in $\phi_q(n)$ and $p$ is in $Z\cap \operatorname{Spec}(\mathcal{O}_{X,q})$. It follows that $p$ is in $\phi_q (n-1)$, and so, $p$ is in $\phi (n-1)$.
\end{proof}

\begin{reminder}\label{rmk:cm_excellent}
    We recall \cite[Definition 1.2]{Cesnavicius:2021}. A locally Noetherian scheme $X$ is \textbf{CM-quasi-excellent} if the following are satisfied:
    \begin{enumerate}
        \item the formal fibers of the local rings of $X$ are Cohen-Macaulay 
        \item every integral closed subscheme contains a nonempty Cohen-Macaulay open subscheme. 
    \end{enumerate}
    Additionally, if $X$ is universally catenary, then we say it is \textbf{CM-excellent}.
\end{reminder} 

\begin{lemma}\label{lem:cm_excellent_covering}
    Let $X$ be a locally Noetherian scheme. The following are equivalent:
    \begin{enumerate}
        \item $X$ is CM-excellent
        \item there exists an open covering of $X$ all of whose members are CM-excellent
        \item for every affine open $\operatorname{Spec}(R)$ of $X$, the ring $R$ is CM-excellent
        \item there exists an affine open covering of $X$ of spectrums of  CM-excellent rings.
    \end{enumerate}
    There exists a similar statement for the case of CM-quasi-excellence.
\end{lemma}

\begin{proof}
    The property of being universally catenary for a locally Noetherian scheme is equivalent to such conditions in the desired claim, see \cite[\href{https://stacks.math.columbia.edu/tag/02J9}{Tag 02J9}]{StacksProject}. Moreover, the property that formal fibers of local rings are Cohen-Macaulay is local, and so the conditions of the desired claim are equivalent in this case. We are left to verify that $(2)$ of \Cref{rmk:cm_excellent} is equivalent to the conditions of the claim. It is evident that $(1)\implies (2)$ in this setting.

    Suppose there exists an open covering of $X$, say $\cup_i U_i$, all of whose members satisfy $(2)$ of \Cref{rmk:cm_excellent}. We want to show that every open affine subscheme of $X$ satisfies $(2)$ of \Cref{rmk:cm_excellent}. Let $j \colon V \to X$ be an open immersion from an affine scheme. Let $Z$ be a closed integral subscheme of $V$ with generic point $\xi$. We can cover by open subscheme (in $Z$) of the form $Z \cap U_i$. The only such $Z\cap U_i$ we are concerned are those that are nonempty, so without of generality, assume this is the case for each $i$. We know that $Z\cap U_i$ corresponds to an integral closed subscheme of $U_i$, which is CM-excellent, so there exists a nonempty Cohen-Macaulay open subscheme of $Z \cap U_i$, say $W_i$. Taking the union $W_i$ gives us a nonempty Cohen-Macaulay subscheme of $Z$, yielding the implication $(2)\implies (3)$.
    
    The case where $(3)\implies (4)$ follows from the fact that a scheme admits an affine open cover, and $(4)\implies (1)$ can be inferred by a similar argument above.
\end{proof}

\begin{example}\label{ex:dualizing_complex_implies_cm_excellent_finite_krull_dim}
    Let $X$ be a Noetherian scheme such that $D_{\operatorname{qc}}(X)$ admits a dualizing complex. Then \cite[\href{https://stacks.math.columbia.edu/tag/0A86}{Tag 0A86}]{StacksProject} tells us $D_{\operatorname{qc}}(U)$ admits a dualizing complex for each affine open in $X$. It follows from \cite[Corollary 1.4]{Kawasaki:2002} that any such $U$ must be CM-excellent of finite Krull dimension. Hence, by \Cref{lem:cm_excellent_covering}, $X$ must be CM-excellent of finite Krull dimension.
\end{example}

\begin{lemma}\label{lem:wc_pointer}
    Let $X$ be a Noetherian scheme. Suppose $Z$ is a closed subset of $X$. If $\phi$ is a Thomason filtration on $X$ which has an associated $t$-structure on $D_{\operatorname{qc}}(X)$ which restricts to $D^b_{\operatorname{coh},Z} (X)$, then $\phi$ is weak Cousin across $Z$.
\end{lemma}

\begin{proof}
    We know from \Cref{cor:descent_for_thomason_on_bounded} and \Cref{lem:weak_cousin_iff_local} that this is an affine local problem. It follows from \cite[Theorem 4.4]{AlonsoTarrio/JeremiasLopez/Saorin:2010}, for the affine case, that $\phi$ must be weak Cousin across $Z$ as this is a problem about support. In fact, the argument for Corollary 4.5 loc.\ cit.\ is readily adaptable in this setting.
\end{proof}

\begin{lemma}[Takahashi]\label{lem:relative_Takahashi}
    Let $R$ be a CM-excellent ring of finite Krull dimension. Suppose $Z$ is a closed subset of $X:=\operatorname{Spec}(R)$. A Thomason filtration $\phi$ on $X$ supported on $Z$ has an associated $t$-structure on $D_{\operatorname{qc},Z}(X)$ which restricts to $D^b_{\operatorname{coh},Z} (R)$ if, and only if, $\phi$ is weak Cousin across $Z$.
\end{lemma}

\begin{proof}
    This argument is strictly a question about support, and hence, it follows by a similar strategy as in the case of \cite[Theorem 5.5]{Takahashi:2023} (which is an adaptation of \cite[Theorem 6.9]{AlonsoTarrio/JeremiasLopez/Saorin:2010} in the presence of dualizing complexes).
\end{proof}

\begin{proof}
    [Proof of \Cref{prop:weak_cousin_cm_excellent_stalks}]
    That the converse direction holds follows from  \Cref{lem:wc_pointer}. 
    We check the forward direction, and this is where our constraint on finite Krull dimension comes to play. Let $\phi$ be a Thomsason filtration on $X$ supported on $Z$ which is weak Cousin across $Z$. Suppose $j\colon U \to X$ is an open immersion from an affine scheme with $Z\cap U\not=\emptyset$. Note that $U$ is CM-excellent, see \Cref{lem:cm_excellent_covering}. By \Cref{lem:weak_cousin_iff_local}, $\phi$ restricted to $U$ will satisfy the weak Cousin across $Z\cap U$. This ensures that the associated $t$-structure of $\phi\cap U$ on $D_{\operatorname{qc},Z\cap U}(U)$ will restrict to a $t$-structure on $D^b_{\operatorname{coh},Z\cap U}(U)$ as $U$ has finite Krull dimension, see \Cref{lem:relative_Takahashi}. Choose an object $E$ in $D^b_{\operatorname{coh},Z}(X)$. Consider the truncation triangle $E$ with respect to associated $t$-structure for $\phi$:
    \begin{displaymath}
        A \to E \to B \to A[1].
    \end{displaymath}
    The observation above tells us that both $A,B$ have bounded and coherent cohomology locally on each affine open subscheme of $X$, both of which are supported on $Z$. It follows from \Cref{lem:local_to_global_aisle} that $A$ and $B$ belongs to $D^b_{\operatorname{coh},Z}(X)$, ensuring that $E$ must as well. This completes the proof.
\end{proof}

\begin{proof}
    [Proof of \Cref{thm:weak_cousin_cm_excellent_stalks}]
    First, we show for any $\otimes$-aisle $\mathcal{A}$ on $D^b_{\operatorname{coh},Z}(X)$, there is a unique Thomason filtration $\phi_{\mathcal{A}}$ on $X$ such that $\mathcal{A} = D^b_{\operatorname{coh},Z}(X)\cap\mathcal{U}_\phi$. By \cite[Theorem 3.10]{Hrbek/Lank/LeGros/Pavon:2026}, $\overline{\langle \mathcal{A} \rangle}^{(-\infty,0]}_{\otimes}$ is compactly generated on $D_{\operatorname{qc}}(X)$. Set $\phi_{\mathcal{A}}$ to be the associated Thomason filtration of $\overline{\langle \mathcal{A} \rangle}^{(-\infty,0]}_{\otimes}$ on $D_{\operatorname{qc}}(X)$ (see \Cref{rem:classification}). By construction, we know that $\phi_{\mathcal{A}}(n)\subseteq Z$ for all $n\in \mathbb{Z}$ and that $\overline{\langle \mathcal{A} \rangle}^{(-\infty,0]}_{\otimes} = \mathcal{U}_{\phi_{\mathcal{A}}}$. It is easy to see that $\mathcal{A} \subseteq \mathcal{U}_{\phi_{\mathcal{A}}} \cap D^b_{\operatorname{coh},Z}(X)$. We check the reverse inclusion. Let $E\in \mathcal{U}_{\phi_{\mathcal{A}}} \cap D^b_{\operatorname{coh},Z}(X)$. Consider the truncation triangle of $E$ with respect to $\mathcal{A}$,
    \begin{displaymath}
        \tau^{\leq 0} E \to E \to \tau^{\geq 1} E \to (\tau^{\leq 0} E)[1].
    \end{displaymath}
    This distinguished triangle corresponds to the truncation triangle for $E$ with respect to $\mathcal{U}_{\phi_{\mathcal{A}}}$ (see e.g.\ proof of \Cref{lem:dbcoh_tensor_via_perfect_tensor}). Since $E\in \mathcal{U}_{\phi_{\mathcal{A}}}$, it follows that the morphism $\tau^{\leq 0} E \to E$ is an isomorphism. Hence, $\mathcal{A} = \mathcal{U}_{\phi_{\mathcal{A}}} \cap D^b_{\operatorname{coh},Z}(X)$.

    Now, for any pair of distinct $\otimes$-aisles $\mathcal{A}$ and $\mathcal{B}$ on $D^b_{\operatorname{coh},Z}(X)$, one has $\mathcal{U}_{\phi_{\mathcal{A}}} \not= \mathcal{U}_{\phi_{\mathcal{B}}}$. Indeed, if $\mathcal{U}_{\phi_{\mathcal{A}}} \not= \mathcal{U}_{\phi_{\mathcal{B}}}$, then 
    \begin{displaymath}
        \mathcal{A} = \mathcal{U}_{\phi_{\mathcal{A}}} \cap D^b_{\operatorname{coh},Z}(X) = \mathcal{U}_{\phi_{\mathcal{B}}} \cap D^b_{\operatorname{coh},Z}(X) = \mathcal{B},
    \end{displaymath}
    which is absurd. Thus, the assignment $\mathcal{A} \to \phi_{\mathcal{A}}$ is injective. By \Cref{prop:weak_cousin_cm_excellent_stalks}, $\phi_{\mathcal{A}}$ satisfies the weak Cousin condition across $Z$. Conversely, given any Thomason filtration $\phi$ on $X$ satisfying the weak Cousin across $Z$, \Cref{prop:weak_cousin_cm_excellent_stalks} yields that $\mathcal{U}_\phi$ is an $\otimes$-aisle which restricts to $D^b_{\operatorname{coh},Z}$. Consequently, we have shown that the classification in \Cref{rem:classification} induces the desired bijection.
\end{proof}

\section{Restriction to perfects}
\label{sec:results2}

Our primary goal in the remainder of this section is to establish \Cref{thm:appearance_for_restrictability_via_filtration}. We start with a few lemmas. 

\begin{lemma}\label{lem:A2}
    Let $(R,\mathfrak{m})$ be a Noetherian local ring. Consider an object $F$ in $D_{\operatorname{qc}}^{\geq n^\prime}(R) \cap D_{\operatorname{qc}}^{\leq n}(R)$ whose support has positive dimension. Suppose the dimension of $\operatorname{supp}(F)$ and $\operatorname{supp}(H^n(F))$ coincide. Then $H^{n+r}_{\mathfrak{m}}(F)$ is not a finitely generated over $R$ where $r:= \operatorname{dim}\operatorname{supp}(F)$.
\end{lemma}

\begin{proof}
    Without loss of generality we can assume $n^\prime = 0$ as
    $H^i_{\mathfrak{m}}(F)$ = $H^{i+j}_{\mathfrak{m}}(F[-j])$ for all $i,j\in \mathbb{Z}$. Recall that $H^i_{\mathfrak{m}}(F)$ is the $i$-th-cohomology of the following complex:
    \begin{displaymath}
        F \overset{\mathbf{L}}{\otimes} \big( R \to \prod_{i_0}R_{f_{i_0}} \to \prod_{i_0 < i_1}R_{f_{i_0}f_{i_1}} \to R_{f_1\cdots f_r} \big)
    \end{displaymath}
    where $f_1,\ldots,f_r$ is a minimal set of generators for $\mathfrak{m}$, see \cite[\href{https://stacks.math.columbia.edu/tag/0A6R}{Tag 0A6R}]{StacksProject}. Denote by $A$ the complex above that we are tensoring $F$ with. Note that $A$ is $K$-flat. It is possible to calculate the cohomology of the tensor product using the spectral sequence coming from the double complex. As the double complex is concentrated in the first quadrant, the spectral sequence converges. Just to fix our convention, we make explicit the zeroth page of the spectral sequence: $E_0^{\ast,j}= F \otimes A^j$ for all $j$ where $A^j$ is the $j$-th component of $A$. The differentials in the first page correspond to the vertical differentials of the double complex. It is easy to see that $E_2^{i,j}=H^j_{\mathfrak{m}}(H^i(F))$. Now, as $H^i_{\mathfrak{m}}(M) = 0$ for all $R$-module $M$ and $i > \operatorname{supp}(M)$, we get that $E_2^{i,j}=0$ if $j > r$ or $i > n$, and so, it follows that:
    \begin{displaymath}
        H^j_{\mathfrak{m}}(H^i(F))=E_2^{n,r} = E_{\infty}^{n,r}=H^{n+r}_{\mathfrak{m}}(F).
    \end{displaymath}
    Thus, $H^j_{\mathfrak{m}}(H^i(F))$ is not a finitely generated $R$-module, see \cite[Proposition A.2]{Smith:2022}
\end{proof}

\begin{lemma}\label{lem:restrict_to_non_regular_local_ring}
    Let $(R,\mathfrak{m},k)$ be a non-regular Noetherian local ring, $Z$ be a closed subset of $\operatorname{Spec}(R)$, and $\phi$ be a Thomason filtration on $\operatorname{Spec}(R)$ which is not constant on the empty set. If the t-structure associated to $\phi$ restricts to $\operatorname{Perf}_Z(R)$, then $\phi$ contains the constant filtration on $\{\mathfrak{m}\}$.
\end{lemma}

\begin{proof}
    For the sake of contradiction, let us suppose that $\phi$ does not contain the constant filtration on $\mathfrak{m}$. Since $\phi$ is not the constant filtration on the the empty set, and does not contain the constant filtration on $\mathfrak{m}$, there is a largest integer $n$ such that $\phi(n) \ne \emptyset$ and $\phi(n+1) = \emptyset$. Since $\phi(n)$ is specialization closed and non-empty, it follows that $\mathfrak{m}\in \phi(n)$. Denote by $K(\mathfrak{m})$ the Koszul complex associated to some sequence of elements generating $\mathfrak{m}$ in $R$. It follows from \Cref{cor:filtrations_agreeing_after_some_point_have_same_truncation} (and \Cref{ex:standard_t_structure_filtration}) that one has a string of isomorphisms:
    \begin{displaymath}
        \begin{aligned}
            \tau^{\geq 0}_\phi (K(\mathfrak{m})[-n]) 
            &\cong \tau^{\geq n}_{\operatorname{st}}(K(\mathfrak{m})[-n]) \\&\cong \tau^{\geq 0}_{\operatorname{st}}(K(\mathfrak{m}))[-n]
            \\& \cong k[-n].
        \end{aligned}
    \end{displaymath}
    However, this tells us that the residue field $k$ is in $\operatorname{Perf}(R)$, and yet $R$ is not regular, giving us a contradiction.
\end{proof}

\begin{proposition}\label{prop:reduce_to_almost_curves}
    Let $(R,\mathfrak{m},k)$ be a Noetherian local ring, $\mathfrak{q}\rightsquigarrow \mathfrak{m}$ be a direct generalization, 
    and $\phi$ be a Thomason filtration on $\operatorname{Spec}(R)$. Then $\phi(i) = \{\mathfrak{m}\}$ for at most finitely many integers $i$ if the associated to a $t$-structure of $\phi$ restricts to $\operatorname{Perf}_{\overline{\{ \mathfrak{q}\} }}(R)$. Additionally, if $\phi$ is not a constant filtration, then $R$ is regular.
\end{proposition}

\begin{proof}
    We can impose $\phi$ be nonempty. Consider the Koszul complex $K(\mathfrak{q})$ associated to $\mathfrak{q}$ in $R$. Note that $K(\mathfrak{q})$ is cohomologically concentrated some the interval $[-n,0]$. Our hypothesis that the associated to a $t$-structure of $\phi$ restricts to $\operatorname{Perf}_{\overline{\{ \mathfrak{q}\} }}(R)$ implies $\tau_\phi^{\leq 0} ( K(\mathfrak{q}) [m] )$ belongs to $\operatorname{Perf}(R)$ for all $m$. It follows from \Cref{lem:A2}, coupled with \Cref{cor:cohomology_for_filtrations_agreeing_on_strip}, that we cannot have $\phi(-n-m) = \phi(-m + 1)= \{\mathfrak{m}\}$ for all integers $m$. This shows the first claim.

    We check the second claim. Assume the contrary, i.e.\ $R$ is not a regular ring. Then \Cref{lem:restrict_to_non_regular_local_ring} tells us that $\phi$ contains the constant filtration on $\{\mathfrak{m}\}$, which contradicts the first claim.
\end{proof}

\begin{lemma}\label{lem:behavior_for_supports_intersection_union}
    Let $Z$ be a connected closed subset of a Noetherian scheme $X$. Suppose that $\phi$ is a Thomason filtration on $X$ whose associated $t$-structure restricts to $\operatorname{Perf}_Z(X)$. Consider a direct generalization $q\rightsquigarrow p$ in $Z$.
    \begin{enumerate}
        \item If $q$ is not in $\bigcap_{n\in \mathbb{Z}} (\phi \cap Z)(n)$, then $p$ is not in $\bigcap_{n\in \mathbb{Z}} (\phi\cap Z)(n)$. 
        \item If $p$ is in $\bigcup_{n\in \mathbb{Z}} (\phi \cap Z)(n)$, then $q$ is in $\bigcup_{n\in \mathbb{Z}} (\phi \cap Z)(n)$. 
    \end{enumerate}
\end{lemma}

\begin{proof}
    We check the fist claim. Consider the Thomason filtration $\psi:=(\phi\cap \overline{\{q\}})$ on $X$. Note that $q$ not belonging to $\bigcap_{n\in \mathbb{Z}} (\phi\cap Z)(n)$ ensures $\psi$ is not constant. It follows from \Cref{lem:all_associated_t_structures_restrict}, coupled with our hypothesis on $\phi$, that $\psi$ has an associated $t$-structure on $D_{\operatorname{qc}}(X)$ which restricts to $\operatorname{Perf}_{\overline{\{q\}}}(X)$. Set $Z_p:=\overline{\{q\}}\cap \operatorname{Spec}(\mathcal{O}_{X,p})$. Then \Cref{prop:reduction_to_affines_stalks} tells us that $\psi_p$ has an associated $t$-structure on $D_{\operatorname{qc},Z_p}(\mathcal{O}_{X,p})$ which restricts to $\operatorname{Perf}_{Z_p} (\mathcal{O}_{X,p})$. Hence, by \Cref{prop:reduce_to_almost_curves}, we have that $p$ does not belong to $\psi_p(n)$ for $n\gg 0$. This tells us $p$ cannot belong to $(\phi\cap Z)(n)$ for $n\gg 0$ as desired. 
    
    The second claim can be checked by arguing in a similar by proving its contrapositive. Specifically, in such a case, $\phi\cap \overline{\{ q\}}$ is not constant from the hypothesis, and one appeals to Lemma ~\ref{lem:restrict_to_non_regular_local_ring} and \Cref{prop:reduce_to_almost_curves} to furnish the claim.
\end{proof}

\begin{lemma}\label{lem:absolute_behavior_union_intersection}
    Let $Z$ be a connected closed subset of a Noetherian scheme $X$. Suppose that $\phi$ is a Thomason filtration on $X$ whose associated $t$-structure restricts to $\operatorname{Perf}_Z(X)$. Then $\bigcup_{n\in \mathbb{Z}} (\phi \cap Z)(n)$ is either $\emptyset$ or $Z$, and similarly for $\bigcap_{n\in \mathbb{Z}} (\phi \cap Z)(n)$.
\end{lemma}

\begin{proof}
    As $Z$ is connected, any pair of irreducible components $Z_{s_1},Z_{s_2}$ enjoy the property that there is a sequence of irreducible components $Z_{v_1}:=Z_s,Z_{v_2},\ldots, Z_{v_k}:=Z_t$ for $Z$ such that $Z_{v_b}\cap Z_{v_{b+1}}\not= \emptyset$. Then the claim follows from \Cref{lem:behavior_for_supports_intersection_union}.
\end{proof}

\begin{proposition}\label{prop:restrictability_implies_regular_locus_containment}
    Let $Z$ be a connected closed subset of a Noetherian scheme $X$. Suppose $\phi$ is a Thomason filtration on $X$ whose associated $t$-structure on $D_{\operatorname{qc}}(X)$ restricts to a $t$-structure on $\operatorname{Perf}_Z(X)$. Then $(\phi\cap Z)(n)$ is contained in the regular locus of $X$ for all $n$ if $\phi$ is not constant on $Z$.
\end{proposition}

\begin{proof}
    We prove this by contradiction. Assume there is an $n$ and $p$ in $(\phi\cap Z)(n)$ which is not contained in the regular locus of $X$. It follows from \Cref{lem:restrict_to_non_regular_local_ring} that $p$ is in $(\phi\cap Z)(t)$ for all $t$ as $(\phi \cap Z)_p$ is not constant on $\operatorname{Spec}(\mathcal{O}_{X,p})$ and has an associated $t$-structure which restricts to $\operatorname{Perf}_{Z\cap \operatorname{Spec}(\mathcal{O}_{X,p})} (\mathcal{O}_{X,p})$ by \Cref{prop:reduction_to_affines_stalks}. Then $\bigcap_{n\in \mathbb{Z}} (\phi \cap Z)(n)$ is nonempty, which ensures by \Cref{lem:absolute_behavior_union_intersection} that $\bigcap_{n\in \mathbb{Z}} (\phi \cap Z)(n) = Z$. However, we see that $\phi$ must be constant on $Z$, which is absurd and completes the proof.
\end{proof}

\begin{corollary}\label{cor:restrictability_implies_weak_cousin_across_and_contained_in_regular_locus}
    Let $Z$ be a connected closed subset of a Noetherian scheme $X$. Suppose $\phi$ is a Thomason filtration on $X$ whose associated $t$-structure on $D_{\operatorname{qc}}(X)$ restricts to a $t$-structure on $\operatorname{Perf}_Z(X)$. Then $\phi$ is weak Cousin across $Z$.
\end{corollary}

\begin{proof}
    The case $\phi$ is constant on $Z$ is obvious, and so, we can impose this is not the case.
    It follows by \Cref{prop:restrictability_implies_regular_locus_containment} that $(\phi\cap Z)(n)$ is contained in the regular locus of $X$ for all $n$. Choose $t$ such that $(\phi\cap Z)(t)$ is nonempty. Let $p$ be in $(\phi \cap  Z)(t)$. Then $\mathcal{O}_{X,p}$ is a regular local ring, and so, $\operatorname{Perf}_{Z\cap \operatorname{Spec}(\mathcal{O}_{X,p})}(\mathcal{O}_{X,p}) = D^b_{\operatorname{coh, }Z\cap \operatorname{Spec}(\mathcal{O}_{X,p})}(\mathcal{O}_{X,p})$. We know by \Cref{prop:reduction_to_affines_stalks} that $\phi \cap \operatorname{Spec}(\mathcal{O}_{X,p})$ has an associated $t$-structure on $D_{\operatorname{qc}}(\mathcal{O}_{X,p})$ which restricts to $\operatorname{Perf}_{Z\cap \operatorname{Spec}(\mathcal{O}_{X,p})}(\mathcal{O}_{X,p})$. However, $\mathcal{O}_{X,p}$ is a regular local ring, and so, $\operatorname{Perf}_{Z\cap \operatorname{Spec}(\mathcal{O}_{X,p})}(\mathcal{O}_{X,p}) = D^b_{\operatorname{coh, }Z\cap \operatorname{Spec}(\mathcal{O}_{X,p})}(\mathcal{O}_{X,p})$. It follows by \Cref{lem:wc_pointer} that $\phi\cap \operatorname{Spec}(\mathcal{O}_{X,p})$ satisfies weak Cousin across $Z\cap \operatorname{Spec}(\mathcal{O}_{X,p})$. Observe that this holds for all $p$ in $\bigcup_{n\in \mathbb{Z}} (\phi\cap Z)(n)$, and so, \Cref{lem:weak_cousin_iff_local} ensures $\phi$ must be weak Cousin across $Z$.
\end{proof}

\begin{lemma}\label{lem:splitting_for_truncations}
    Let $Z$ be a closed subset of a Noetherian scheme. Denote the connected components of $Z$ by $Z_t$. Suppose $\phi$ is a Thomason filtration on $X$. Then $\tau^{\leq 0}_\phi (E)$ is isomorphic to $\bigoplus_t \tau^{\leq 0}_{\phi\cap Z_t} E$ and $\tau^{\geq 1}_\phi (E)$ is isomorphic to $\bigoplus_t \tau^{\geq 1}_{\phi\cap Z_t} E$ if $E$ is an object of $D^b_{\operatorname{coh},Z}(X)$.
\end{lemma}

\begin{proof}
    We know that the connective and coconnective truncation functors commute with finite coproducts. It follows by \cite[Lemma 3.9]{Huybrechts:2006} that there are objects $E_t$ in $D^b_{\operatorname{coh},Z_t}(X)$ such that $E$ is isomorphic to $\bigoplus_t E_t$. There is a string of isomorphisms:
    \begin{displaymath}
        \begin{aligned}
            \tau^{\leq 0}_\phi (E) &\cong \tau^{\leq 0}_\phi (\bigoplus_t E_t) \\&\cong \bigoplus_t \tau^{\leq 0}_\phi E_t \\&= \bigoplus_t \tau^{\leq 0}_{\phi\cap Z_t} E_t && (\textrm{by \Cref{prop:truncations_preserve_support}}).
        \end{aligned}
    \end{displaymath} 
    It remains to show $\tau^{\leq 0}_{\phi\cap Z_t} E$ is isomorphic to $\tau^{\leq 0}_{\phi\cap Z_t} E_t$ for each $t$. We have another string of isomorphisms:
    \begin{displaymath}
        \begin{aligned}
            \tau^{\leq 0}_{\phi\cap Z_t} (E) &\cong \tau^{\leq 0}_{\phi\cap Z_t} (\bigoplus_s E_s) \\&\cong \bigoplus_s \tau^{\leq 0}_{\phi\cap Z_t} E_s \\&= \bigoplus_s \tau^{\leq 0}_{\phi\cap Z_t\cap Z_s} E_s && (\textrm{by \Cref{prop:truncations_preserve_support}}).
        \end{aligned}
    \end{displaymath}
    Note that $Z_s\cap Z_t=\emptyset$ for $s\not=t$, and so, $\tau^{\leq 0}_{\phi\cap Z_t\cap Z_s} E_s=0$ in such cases. Hence, the desired claim for connective truncations has been shown, and a similar argument gives the claim for coconnective truncations. 
\end{proof}

\begin{lemma}\label{lem:restriction_iff_to_connected_components}
    Let $Z$ be a closed subset of a Noetherian scheme with connected components $Z_t$. A Thomason filtration $\phi$ on $X$ has associated $t$-structure on $D_{\operatorname{qc}}(X)$ that restricts to $\operatorname{Perf}_Z (X)$ if, and only if, it restricts to $\operatorname{Perf}_{Z_t} (X)$ for each $t$.
\end{lemma}

\begin{proof}
    This follows from \Cref{cor:automatic_restriction_to_smaller_Thomason} and \Cref{lem:splitting_for_truncations}.
\end{proof}

\begin{proof}
    [Proof of \Cref{thm:appearance_for_restrictability_via_filtration}]
    It follows by \Cref{lem:restriction_iff_to_connected_components} that we can impose $Z$ be connected. 

    Suppose $Z$ is contained in the regular locus of $X$. Then $\operatorname{Perf}_Z (X) = D^b_{\operatorname{coh},Z}(X)$. It follows from \Cref{cor:lifitng_tensor_structure_from_perfect_categories} that any tensor $t$-structure on $\operatorname{Perf}_Z (X)$ is the restriction of a $t$-structure on $D_{\operatorname{qc},Z}(X)$ that is associated to a Thomason filtration on $X$. This case is determined by the Thomason filtration in \Cref{ex:standard_t_structure_filtration}, which shows the forward direction.

    Now we check the converse direction. Suppose there is a Thomason filtration on $X$ which is nonconstant on $Z$ and whose associated $t$-structure on $D_{\operatorname{qc}}(X)$ restricts to $\operatorname{Perf}_Z (X)$. Then \Cref{prop:restrictability_implies_regular_locus_containment} ensures that $Z$ must be in regular locus, which completes the proof.
\end{proof}

\begin{example}\label{ex:nagata_singular}
    \hfill
    \begin{enumerate}
        \item There is an example of a regular Noetherian ring $R$ of infinite Krull dimension due to Nagata (see \cite[Example 1]{Nagata:1962}). Let $S$ be any Noetherian ring which is not regular (e.g.\ $k[ x,y,z ] / (x^2 y^7 z^{11})$). Then $S\oplus R$ is a Noetherian ring of infinite Krull dimension which cannot be regular.
        \item We construct a more `complicated' instance. See \cite[Proposition 1]{Hochster:1973} for details on the following construction. Let $k$ be a field. Suppose $R$ is a geometrically integral\footnote{See \cite[\href{https://stacks.math.columbia.edu/tag/05DW}{Tag 05DW}]{StacksProject}.} $k$-algebra which is of finite type over $k$ and not regular (e.g.\ $\mathbb{C}[w,x,y,z]/(wz-xy)$). Set $R_1:=R$ and $R_n = k[x_1,\ldots,x_n]$ for all $n\geq 2$. Consider the maximal ideals $\mathfrak{p}_1:=(x,y,z,w)$ of $R_1$ and $\mathfrak{p}_n = (x_1,\ldots,x_n)$ of $R_n$ for $n\geq 2$. Denote by $R$ for the tensor product of all the $R_n$'s over $k$. The set $S:=R\setminus (\cup_{n\geq 1} \mathfrak{p}_n R)$ is multiplicatively closed in $R$. Let $R^\prime$ be the localization of $R$ at $S$. Then $R^\prime$ is a Noetherian integral domain whose maximal ideals correspond to $\mathfrak{p}_n R^\prime$. Our choices tells us, by \cite[Proposition 1]{Hochster:1973}, that $R^\prime$ must have infinite Krull dimension and cannot be regular (as localization at $\mathfrak{p}_1 R^\prime$ is not regular). 
    \end{enumerate}
\end{example}

\begin{proof}
    [Proof of \Cref{cor:classification_for_relative_perf}]
    We can reduce, by \Cref{lem:restriction_iff_to_connected_components}, to the case where $Z$ is connected. Then the claim follows by \Cref{prop:weak_cousin_cm_excellent_stalks}
    and \Cref{thm:appearance_for_restrictability_via_filtration}.
\end{proof}

\bibliographystyle{alpha}
\bibliography{mainbib}

\end{document}